\documentclass[12pt]{article}
\usepackage[latin1]{inputenc}
\usepackage[T1]{fontenc}
\usepackage{a4wide}
\usepackage{fontenc}
\usepackage{color}
\usepackage{amsmath}
\usepackage{graphicx}
\usepackage{array}
\usepackage{amsfonts}
\usepackage{amssymb}
\usepackage{euscript}
\usepackage{delarray}
\usepackage{latexsym}
\usepackage{enumerate}
\usepackage{float}

\title{Regularity dependence of the rate of convergence of the learning curve  for Gaussian process regression.}
\author{Loic Le Gratiet $^\dag$  $^\ddag$, Josselin Garnier $^*$  \\ $^\dag$ Universit\'e Paris Diderot 75205 Paris Cedex 13 \\ $^\ddag$ CEA, DAM, DIF, F-91297 Arpajon, France  \\
$^*$ Laboratoire de Probabilites et Modeles Aleatoires \& \\ Laboratoire Jacques-Louis Lions,\\ Universite Paris Diderot, 75205 Paris Cedex 13, France}

\newtheorem{theo}{Theorem}
\newtheorem{example}{Example}
\newtheorem{prop}{Proposition}

\begin{document}

\maketitle

\section{Abstract}

This paper deals with   the learning curve in a Gaussian process regression framework. The learning curve describes the   generalization error of the Gaussian process used for the regression.  The main result is the  proof of a theorem giving the  generalization error for a large class of correlation kernels and for  any dimension when the number of observations is large. From this result, we can deduce the asymptotic behavior of  the generalization error when the observation error is small. The presented proof generalizes previous ones   that were limited to  special kernels or to small dimensions (one or two).
The result can be used to build an optimal strategy for resources allocation. This strategy is applied successfully to a nuclear safety problem.

\paragraph{Keywords:}  Gaussian process regression, asymptotic mean squared error, learning curves, generalization error, convergence rate.

\section{Introduction}

Gaussian process regression is a useful tool to approximate an objective function given some of its  observations \cite{Las94}.
 It has originally been used in  geostatistics to interpolate  a random field at  unobserved locations \cite{Wa03}, \cite{Ber01} and \cite{Gnei10}, it has been developed in many areas such as  environmental and atmospheric sciences.

 This method has become very  popular during the last decades  to build surrogate models from noise-free observations. For example, it is widely used in the field of ``computer experiments'' to build models which surrogate an expensive computer  code \cite{SACKS89}. Then, through the fast approximation of the computer code, uncertainty  quantification and sensitivity analysis can be performed with a low computational cost.

Nonetheless, for many realistic cases, we do not have direct access to the function to be approximated  but only to noisy versions of it. For example, if the objective function is the result of an experiment, the available responses  can be tainted by  measurement noise. In that case, we can reduce the noise of the observations by repeating the experiments at the same locations. Another example is   Monte-Carlo based simulators - also called stochastic simulators - which use  Monte-Carlo or Monte-Carlo Markov Chain methods to solve a system of differential equations through its probabilistic interpretation. For such simulators, the noise level can be tuned by the number of Monte-Carlo particles used in the procedure.

Gaussian process regression can be easily adapted to the case of noisy observations.
The purpose of this paper is to minimize the generalization error - defined as the averaged mean squared error - of the Gaussian process regression with noisy observations and a given budget. The budget is defined as the number of experiments including the repetitions.
As seen in the previous paragraph, in many cases  the noise variance is inversely proportional to the number of repetitions. Therefore, if the total budget is given, a trade off between the number and the accuracy of the observations has to be made. 

Many authors were interested in obtaining learning curves describing the generalization error as a function of the training set size \cite{R06}. The problem has been addressed in the statistical and numerical analysis areas. For an overview, the reader is referred to \cite{RiY00} for a numerical analysis point of view and to \cite{R06} for a statistical one. In particular, in the numerical analysis literature, the authors are interested in numerical differentiation of functions from noisy data (see \cite{Ri96} and \cite{Bo03}). They have found very interesting results for kernels satisfying the Sacks-Ylvisaker conditions of order $r$ \cite{SY89} but only valid for 1-D or 2-D functions. 

In the statistical literature \cite{SH02} give accurate approximations to the learning curve and \cite{OV99} and \cite{WV00} give upper and lower bounds on it. Their approximations give the asymptotic value of the learning curve (for a very large number of observations).  They are based on the Woodbury-Sherman-Morrison matrix inversion lemma \cite{Har97} which holds in finite-dimensional cases which correspond to degenerate covariance kernels in our context. Nonetheless, classical kernels used in Gaussian process regression are non-degenerate and we hence are in an infinite-dimensional case and the Woodbury-Sherman-Morrison formula cannot be used directly. Another proof for degenerate kernels can be found in  \cite{VP09}.

The main result of this paper is a theorem giving the   value of the Gaussian process regression mean squared error for a large training set size when the  observation noise variance is proportional to  the number of observations. This   value is given as a function of the eigenvalues and  eigenfunctions   of the  covariance kernel.  From this theorem, we can deduce an approximation of the learning curve for non-degenerate and degenerate kernels (which generalizes  results in \cite{OV99}, \cite{SH02} and \cite{VP09}) and for any dimension (which generalizes results in  \cite{RiY00}, \cite{Ri96} and \cite{Bo03}). Finally, from this approximation we can deduce the rate of convergence of the Best Linear Unbiased Predictor (BLUP) in a  Gaussian process regression framework.

The rate of convergence of the BLUP is of practical interest since it provides a powerful tool for decision support. Indeed, from an initial experimental design set, it can  predict the additional computational budget necessary  to reach   a given desired accuracy when the  observation noise variance is homogeneous in space.  Finally, we propose in this paper a theorem giving the best  resource allocation when the noise variance is heterogeneous in space.

The paper is organized as follow. First we present the considered Gaussian process regression model with noisy observations. Second, we present the main result of the paper which is the theorem giving the mean squared error  of the considered model for a large training size. Third, we study the rate of convergence of the generalization error when the noise variance decreases. Academic examples are presented   to compare the theoretical convergences given by the theorem and  numerically observed convergences. Then, we address the problem of the optimal budget allocation.  Finally, an industrial application to  the safety assessment of a nuclear system containing fissile materials is considered. This real case emphasizes the effectiveness of the theoretical  rate of convergence of the BLUP since it predicts a very good approximation of the budget needed  to reach a prescribed precision.

\section{Gaussian process regression}

Let us suppose that we want to approximate an objective function   $x \in \mathbb{R}^d \rightarrow f(x) \in \mathbb{R}$ from noisy observations of it at points $(x_i)_{i=1,\dots,n}$ with $x_i \in \mathbb{R}^d$. The points of the experimental design set $(x_i)_{i=1,\dots,n}$ are supposed to be sampled from the probability measure $\mu$ over $\mathbb{R}^d$.  $\mu$ is called  the design measure, it can have either a compact support (for a bounded input parameter space  domain) or  unbounded support (for   unbounded input parameter space). We hence  have $n$ observations of the form $z_i = f(x_i) + \varepsilon(x_i)$ and we consider that $(\varepsilon(x_i))_{i=1,\dots,n}$ are  independently     sampled from the   Gaussian distribution with mean zero and variance $n\tau(x_i)$:
\begin{equation}
\varepsilon(x) \sim \mathcal{N}(0,n\tau(x))
\end{equation}
Note that the number of observations and the  observation noise variance are both controlled by $n$. It means that if we increase the number $n$ of observations, we automatically increase the uncertainty on the observations.  An observation noise variance   proportional to $n$ is natural in the framework of experiments with repetitions or stochastic simulators. Indeed,  for a   fixed number of experiments (or simulations), the user can decide to perform  them in few points with many repetitions (in that case the noise variance will be low) or to perform them in many  points with few repetitions (in that case the noise variance will be large). We introduce in Example \ref{ex1} the framework of repeated experiments. We note that the framework is the same as the one of stochastic simulators and it is the one considered in Sections \ref{allocation} and \ref{application}.

\begin{example} [Gaussian process regression with repeated experiments]\label{ex1}

Let us consider  that we want to approximate the  function   $x \in \mathbb{R}^d \rightarrow f(x) \in \mathbb{R}$ from noisy observations  at points $(x_i)_{i=1,\dots,n}$   sampled from the design measure $\mu$ and  with $s$ replications at each point. We hence  have $ns$ data  of the form $z_{i,j} = f(x_i) + \varepsilon_j(x_i)$ and we consider that $(\varepsilon_j(x_i))_{\substack{i=1,\dots,n \\ j =1,\dots,s}}$ are  independently  distributed from a   Gaussian distribution with mean zero and variance $\sigma_\varepsilon^2(x_i)$. Then, denoting the vector of observed values  by $z^{n} = (z_i^n)_{i=1,\dots,n} = (\sum_{j=1}^s z_{i,j}/s)_{i=1,\dots,n}$, the variance of an observation $z^{n}_i$ is $ \sigma_\varepsilon^2(x_i)/s$. Thus, if we  consider  a fixed budget $T=ns$, we have $ \sigma_\varepsilon^2(x_i)/s = n \tau(x_i)$ with $\tau(x_i) = \sigma_\varepsilon^2(x_i)/T$ and the observation noise variance is   proportional to $n$.

In Section \ref{convergence} we give the value of the generalization error for $n$ large. Then, in Section \ref{rate} we are interested in its    convergence for $n$ large and when $\tau(x)$ tends to zero. Finally, in Section \ref{allocation} we consider the non-uniform allocation $(s_i)_{i=1,\dots,n}$ with $T = \sum_{i=1}^n s_i$ and we   address the question of optimal allocation of   the repetitions $(s_i)_{i=1,\dots,n}$ as a function of the noise level $\sigma_\varepsilon^2(x_i)$ so as to minimize the generalization error.
\end{example}

The main idea of the Gaussian process regression is to suppose that the objective function $f(x)$ is a realization of a Gaussian process $Z(x)$ with a known mean and a known covariance kernel $k(x,x')$. The mean can be considered equal to zero without loss of generality. Then, denoting  by $z^{n} = [f(x_i)+\varepsilon(x_i)]_{1\leq i \leq n}$ the vector of length $n$ containing the noisy observations, we choose as predictor the Best Linear Unbiased Predictor (BLUP) given by the equation:
\begin{equation}\label{BLUP}
\hat{f}(x) = k(x)^T (K+n\Delta)^{-1}z^{n}, \qquad \Delta = \mathrm{diag}[(\tau(x_i))_{i=1,\dots,n}]
\end{equation}
where $k(x) = [k(x,x_i)]_{1\leq i \leq n}$ is the $n$-vector containing the covariances between $Z(x)$ and $Z(x_i), \quad 1\leq i \leq n$ and  $K = [k(x_i,x_j)]_{1\leq i,j \leq n}$ is the $n \times n$-matrix containing the covariances between  $Z(x_i)$ and $Z(x_j), \quad 1\leq i ,j \leq n$. When $\tau(x)$ is independent of $x$, we have $\Delta =\tau I$ with $I$   the $n \times n$ identity matrix. The BLUP minimizes the Mean Squared Error (MSE) which equals:
\begin{equation}\label{MSE_BLUP}
\sigma^2(x) = k(x,x) - k(x)^T (K+n\Delta)^{-1}k(x)
\end{equation}

Indeed, if we consider a Linear Unbiased Predictor (LUP) of the form $a(x)^T z^{n}$, its MSE is given by:
\begin{equation}\label{MSE_LUP}
\mathbb{E}[(Z(x) - a^T(x) Z^{n})^2]  = k(x,x) - 2a(x)^Tk(x) + a(x)^T (K+n\Delta) a(x)
\end{equation}
where $Z^{n} = [Z(x_i)+\varepsilon(x_i)]_{1 \leq i \leq n}$ and $\mathbb{E}$ stands for the expectation with respect to the distribution of the Gaussian process $Z(x)$.
The value of $a(x)$  minimizing (\ref{MSE_LUP}) is  $a_{\mathrm{opt}}(x)^T=k(x)^T (K+n\Delta)^{-1}$. Therefore, the BLUP  given by $a_{\mathrm{opt}}(x)^T z^{n}$ is equal to  (\ref{BLUP}) and  by substituting  $a(x)$ with    $a_{\mathrm{opt}}(x)$ in equation (\ref{MSE_LUP}) we obtain the MSE of the BLUP given by   equation (\ref{MSE_BLUP}).\\

The main result of this paper is the proof of a theorem that gives  the  asymptotic value of $\sigma^2(x)$ when $n \rightarrow +\infty$ and $\Delta = \tau I$. Thanks to this theorem, we can deduce the asymptotic value of the Integrating Mean Squared Error (IMSE) - also called learning curve or generalization error - when  $n \rightarrow +\infty$. The  IMSE is defined by:
\begin{equation}\label{IMSE}
\mathrm{IMSE} = \int_{\mathbb{R}^d} \sigma^2(x) \, d\mu(x)
\end{equation}
where $\mu$ is the design measure of the input space parameters.
The  asymptotic value of the IMSE that we obtain can be viewed as a generalization of previous results (see \cite{R06}, \cite{RiY00}, \cite{Ri96}, \cite{Bo03}, \cite{OV99}, 
\cite{SH02} and \cite{VP09}). It can be used to determine  the  budget required  to reach a prescribed accuracy (see Section \ref{allocation}). Note that the  proof of the theorem holds for a constant  observation noise variance $\tau$.  Nevertheless, to provide optimal resource allocation, it can be important to take into account the heterogeneity of the  observation noise  variance. We give in this paper  under certain restricted conditions (i.e., when $K$ is diagonal) the optimal allocation taking into account the noise heterogeneity. Moreover, we numerically observe that this allocation remains efficient in more general cases although it is not anymore optimal (it remains more efficient than the uniform one).

\section{Convergence of the learning curve for Gaussian process regression}\label{convergence}

This section deals with the convergence of the BLUP when the number of observations is large and   the reduced noise variance does not depend on $x$, i.e. $\tau(x) = \tau$ and $\Delta = \tau I$. The speed of convergence of the BLUP is evaluated through the generalization error - i.e. the IMSE - defined in   (\ref{IMSE}). 
The main theorem of this paper follows:

\begin{theo}\label{theo1}
Let us consider $Z(x)$ a Gaussian process with zero  mean and covariance kernel $k(x,x') \in \mathcal{C}^0(\mathbb{R}^d \times \mathbb{R}^d )$ and $(x_i)_{i=1,\dots,n}$ an experimental design set of $n$ independent random points sampled with the probability measure $\mu$  on $\mathbb{R}^d$. We assume that $\sup_{x \in \mathbb{R}^d} k(x,x) < \infty$. According to Mercer's theorem \cite{Mer09}, we have the following representation of $k(x,x')$:
\begin{equation}
k(x,x') = \sum_{p \geq 0} \lambda_p \phi_p(x) \phi_p(x')
\end{equation}
where $(\phi_p(x) )_p$ is an orthonormal  basis of $L^2_\mu(\mathbb{R}^d)$ (denoting  the set of square integrable functions) consisting of eigenfunctions of $(T_{\mu,k} f)(x)=\int_{\mathbb{R}^d} k(x,x')f(x')d\mu(x')$ and  $\lambda_p$ is the nonnegative sequence of corresponding eigenvalues sorted  in decreasing order. Then, for a non-degenerate kernel - i.e. when $\lambda_p > 0$, $\forall p  > 0$ -     we have the following convergence in probability for the MSE (\ref{MSE_BLUP}) of the BLUP:
\begin{equation}\label{MSE_prob}
\sigma^2(x)  \stackrel{n \rightarrow \infty}{\longrightarrow}  \sum_{p \geq 0} \frac{ \tau \lambda_p}{\tau + \lambda_p}  \phi_p(x)^2
\end{equation}
For degenerate kernels - i.e. when only   a finite number of $\lambda_p$ are not  zero -  the convergence is almost sure.
We note that we have the  convergences   with respect to the distribution of the points $(x_i)_{i=1,\dots,n}$ of the experimental design set.
\end{theo}

The sketch of the proof of Theorem \ref{theo1} is given below. The full proof is given in  Appendix \ref{A_proof}.\\

\emph{Sketch of Proof.} We first prove the theorem for degenerate kernels (see Appendix \ref{A_proof}.1) which was already known in that case. Next we find a lower bound for $\sigma^2(x)$ for non-degenerate kernels. Let us consider the Karhunen-Lo\`eve decomposition of  
$
Z(x) = \sum_{p \geq 0}Z_p \sqrt{\lambda_p} \phi_p(x)
$  where $(Z_p)_p$ is a sequence of independent Gaussian random variables with mean zero and variance one.
If we denote by $a_{\mathrm{opt},i}(x)$, $i=1,\dots,n$, the coefficients of the BLUP associated to $Z(x)$, the Gaussian process regression mean squared error can be written
$
\sigma^2(x) = \sum_{p \geq 0} \lambda_p \left( \phi_p(x) - \sum_{i=1}^{n} a_{\mathrm{opt},i}(x) \phi_p(x_i) \right)^2
$.
Then, for a fixed $\bar{p}$, the following inequality holds:
\begin{equation}\label{ineq_sig_sigLUP}
\sigma^2(x) \geq \sum_{p \leq \bar{p}} \lambda_p \left( \phi_p(x) - \sum_{i=1}^{n} a_{\mathrm{opt},i}(x) \phi_p(x_i) \right)^2 = \sigma^2_{LUP,\bar{p}}(x)
\end{equation}
where, $\sigma^2_{LUP,\bar{p}}(x)$ is the MSE of the Linear Unbiased Predictor (LUP) of coefficients $a_{\mathrm{opt},i}(x)$ associated to the Gaussian process $Z_{\bar{p}}(x) = \sum_{p \leq \bar{p}}Z_p \sqrt{\lambda_p} \phi_p(x)$. Let us consider $\sigma^2_{\bar{p}}(x)$ the MSE of the BLUP of  $Z_{\bar{p}}(x)$, we have  the following inequality: 
\begin{equation}\label{ineq_sigLUP_sigBLUP}
\sigma^2_{LUP,\bar{p}}(x) \geq  \sigma^2_{\bar{p}}(x)
\end{equation}
Since $Z_{\bar{p}}(x)$ has a degenerate kernel, $\forall \bar{p} > 0$, the almost sure convergence (\ref{MSE_prob}) holds for $\sigma^2_{\bar{p}}(x)$. Then, considering inequalities (\ref{ineq_sig_sigLUP}),  the convergence (\ref{MSE_prob})   for $\sigma^2_{\bar{p}}(x)$   and the limit $\bar{p} \rightarrow \infty$, we obtain:
\begin{equation}\label{lowerbound}
\liminf_{n \rightarrow \infty} \sigma^2(x) \geq \sum_{p \geq 0}  \frac{\tau \lambda_p}{\tau +  \lambda_p}  \phi_p(x)^2
\end{equation}

It remains to  find an upper bound for $\sigma^2(x)$. Since $\sigma^2(x)$ is the MSE of the BLUP associated to $Z(x)$, if we consider any other LUP associated to $Z(x)$ , then the corresponding  MSE denoted by $\sigma_{LUP}^2(x)$ satisfies the following inequality:
\begin{displaymath}
\sigma^2(x) \leq \sigma_{LUP}^2(x)
\end{displaymath}
The idea is to find a LUP so that its MSE is a tight upper bound of $\sigma^2(x)$. Let us consider the LUP:
\begin{equation}\label{LUP}
\hat{f}_{LUP}(x) = k(x)^T A z^{n}
\end{equation}
with $A$ the $n \times n$ matrix defined by 
$
A = L^{-1}+\sum_{k=1}^q(-1)^k(L^{-1}M)^kL^{-1}
$
with  $L = n \tau I + \sum_{p < p^*} \lambda_p [\phi_p(x_i) \phi_p(x_j)]_{1 \leq i,j \leq n}$, $M = \sum_{p \geq p^*} \lambda_p [\phi_p(x_i) \phi_p(x_j)]_{1 \leq i,j \leq n}$, $q$ a finite integer and $p^*$ such that $\lambda_{p^*} < \tau$. The choice of this LUP is motivated by the fact that the matrix $A$ is an approximation of the inverse of the matrix $(n \tau I + K)$ that is tractable in the following calculations. Remember that the BLUP is $\hat{f}_{\mathrm{BLUP}}(x) = k(x)^T(K+n \tau I)^{-1}z^{n}$. Then, the  MSE of the LUP (\ref{LUP}) is given by:
\begin{displaymath}
\sigma^2_{LUP}(x) = k(x,x) - k(x)^T L^{-1} k(x) - \sum_{i=1}^{2q+1} (-1)^ik(x)^T(L^{-1}M)^iL^{-1}k(x)
\end{displaymath}
Thanks to the Woodbury-Sherman-Morrison formula\footnote{If $B$ is a non-singular $p\times p$ matrix, $C$ a non-singular $m \times m$ matrix and $A$ a $m \times p$ matrix with $m,p < \infty$, then $(B+AC^{-1}A)^{-1} = B^{-1} - B^{-1}A(A^TB^{-1}A+C)^{-1}A^TB^{-1}$.}, the strong law of large numbers and the continuity of the inverse operator in the space of $p$-dimensional invertible matrices, we have the following almost sure convergence:
\begin{displaymath}
k(x)^T L^{-1}k(x) \stackrel{n \rightarrow \infty}{\longrightarrow} \sum_{p <p^*} \frac{ \lambda_p^2}{ \lambda_p + \tau} \phi_p(x)^2 + \frac{1}{\tau} \sum_{p \geq p^*} \lambda_p^2 \phi_p(x)^2
\end{displaymath}
We note that we can use the Woodbury-Sherman-Morrison formula and the strong law of large numbers since $p^*$ is finite and independent of $n$. Then, using the Markov inequality and the equality $\sum_{p \geq 0} \lambda_p \phi_p(x)^2 = k(x,x) < \infty$, we have the following convergence in probability:
\begin{displaymath}
k(x)^T (L^{-1}M)^iL^{-1}k(x) \stackrel{n \rightarrow \infty}{\longrightarrow} \left( \frac{1}{\tau}\right)^{i+1} \sum_{p \geq p^*} \lambda_p^{i+2}\phi_p(x)^2
\end{displaymath}
We highlight that we cannot use the strong law of large numbers here due to the infinite sum in the definition of $M$. Finally, we obtain the following convergence in probability:
\begin{displaymath}
\limsup_{n \rightarrow \infty} \sigma^2(x) \leq \lim_{n \rightarrow \infty}  \sigma^2_{LUP}(x)  =  \sum_{p \geq 0} \left( \lambda_p - \frac{ \lambda_p^2}{\tau +  \lambda_p} \right) \phi_p(x)^2 - \sum_{p \geq p^*} \lambda_p^2 \frac{\left( \frac{ \lambda_p}{\tau} \right)^{2q+1}}{\tau +  \lambda_p} \phi_p(x)^2
\end{displaymath}
By taking the limit $q \rightarrow \infty$ in the right hand side and using the inequality $\lambda_{p^*} < \tau$, we obtain the following upper bound for $\sigma^2(x)$:
\begin{equation}\label{upperbound}
\limsup_{n \rightarrow \infty} \sigma^2(x) \leq \sum_{p \geq 0}  \frac{\tau \lambda_p}{\tau +  \lambda_p}   \phi_p(x)^2
\end{equation}
The result announced in  Theorem \ref{theo1} is deduced from the lower and upper bounds (\ref{lowerbound}) and  (\ref{upperbound}). $\blacksquare$ 

\paragraph{Remark 1} For non-degenerate kernels such that $|| \phi_p(x) ||_{L^{\infty}} < \infty$ uniformly in $p$, the convergence is almost sure. Some kernels such as the one of the  Brownian motion satisfy this property.\\

The following theorem gives the asymptotic value of the learning curve when $n$ is large.

\begin{theo}\label{theo2}
Let us consider $Z(x)$ a Gaussian process with known mean and covariance kernel $k(x,x') \in \mathcal{C}^0(\mathbb{R}^d \times \mathbb{R}^d )$ such that $\sup_{x \in \mathbb{R}^d} k(x,x) < \infty$  and $(x_i)_{i=1,\dots,n}$ an experimental design set of $n$ independent random points sampled with the probability measure $\mu$   on $\mathbb{R}^d$. Then, for a non-degenerate kernel, we have the following convergence in probability:
\begin{equation}\label{conv_IMSE}
\mathrm{IMSE}  \stackrel{n \rightarrow \infty}{\longrightarrow}  \sum_{p \geq 0} \frac{ \tau \lambda_p}{\tau + \lambda_p} 
\end{equation}
For degenerate kernels, the convergence is almost sure.
\end{theo}

\emph{Proof.} From   Theorem \ref{theo1} and  the orthonormal property of the basis $(\phi_p(x))_p$ in $L^2_\mu(\mathbb{R^d})$, the proof of the theorem is straightforward by integration. We note that we can permute the integral and  the limit thanks to  the dominated convergence theorem since $\sigma^2(x) \leq k(x,x)$. $\blacksquare$ 

\paragraph{Remark 2} The obtained limit is identical to the one established in  \cite{R06} and  \cite{VP09} for a degenerate kernel. Furthermore, \cite{OV99} gives accurate upper and lower bounds for the asymptotic behavior of the IMSE for a degenerate kernel too. The originality of the presented result is the  proof giving the asymptotic value of the learning curve for a non-degenerate kernel.  We note that this result is of practical interest since the usual kernels for Gaussian process regression are non-degenerate and we will exhibit dramatic differences between the learning curves of degenerate and non-degenerate kernels.

\begin{prop}\label{theo3}
 Let us denote  $\mathrm{IMSE}_\infty = \lim_{n \rightarrow \infty}\mathrm{IMSE} $. The following inequality holds:
\begin{equation}\label{rate_IMSE}
\frac{1}{2}B_\tau \leq \mathrm{IMSE}_\infty \leq B_\tau
\end{equation}
with $B_\tau = \sum_{p   \mathrm{ \,\, s.t. \,\,} \lambda_p \leq  \tau} \lambda_p + \tau \# \left\{  p \mathrm{\,\, s.t. \,\,} \lambda_p >  \tau  \right\}$.
\end{prop}
\emph{Proof.} The proof is directly deduced from   Theorem \ref{theo2} and the following inequality:
\begin{displaymath}
\frac{1}{2} h_\tau(x) \leq \frac{x}{x + \tau} \leq h_\tau(x)
\end{displaymath}
with:
\begin{displaymath}
h_\tau(x) = \left\{ \begin{array}{ll} x /\tau& x \leq \tau \\ 1 & x > \tau \\ \end{array} \right.
\end{displaymath}
 $\blacksquare$ 
%
%\paragraph{Remark 3} The inequality given in Proposition \ref{theo3} provides a lower bound more precise in our case than the one given in \cite{MW81} -  $\sum_{p \geq s+1} \lambda_p \leq \mathrm{IMSE}_\infty$. This lower bound corresponds to the case of an experimental design set made with the points of a Gaussian quadrature  whereas in our case  it is sampled randomly  from a given measure.

\section{Examples of rates of convergence for the learning curve}\label{rate}

 Proposition \ref{theo3}  shows that the rate of convergence of the generalization error $\mathrm{IMSE}_\infty $   in function of $\tau$  is  equivalent to  the one of  $B_\tau$. In this Section, we analyze the rate of convergence of $\mathrm{IMSE}_\infty $ (or equivalently $B_\tau$) when $\tau$ is small. We note that the presented results can be interpreted as a rate of convergence in function of the number of observations since $\tau$ is the ratio between the noise variance $n\tau$  and the number of observations $n$.

In this section, we consider that the design measure $\mu$ is uniform on $[0,1]^d$.

\paragraph{Example 2  (Degenerate kernels)}
For degenerate kernels we have $ \# \left\{  p \mathrm{\,\, s.t. \,\,} \lambda_p > 0  \right\} < \infty$. Thus, when $\tau \rightarrow 0$, we have:
\begin{displaymath}
\sum_{p   \mathrm{ \,\, s.t. \,\,} \lambda_p < \tau} \lambda_p  = 0
\end{displaymath}
from which:
\begin{equation}\label{T_degenerate}
B_\tau \propto \tau
\end{equation}

Therefore, the IMSE decreases as $\tau$. We find here a classical result about Monte-Carlo convergence which gives that the variance decay is   proportional to the observation noise variance  ($n\tau$) divided by the number of observations ($n$)   whatever the dimension. Nevertheless, for non-degenerate kernels, the number of non-zero eigenvalues is infinite and we are hence in an infinite-dimensional case (contrarily to the degenerate one). We see in the following examples that we do not conserve the usual Monte-Carlo convergence rate in this case which emphasizes  the importance of  Theorem \ref{theo1} dealing with non-degenerate kernels.

\paragraph{Example 3  (The fractional Brownian motion)} 
Let us consider the fractional Brownian kernel with Hurst parameter $H  \in (0,1)$:
\begin{equation}\label{fBm_kernel}
k(x,y) = x^{2H} + y^{2H} - |x-y|^{2H}
\end{equation}
The  associated Gaussian process - called fractional Brownian motion -  is H\"older continuous with exponent $H-\varepsilon$,  $\forall \varepsilon > 0$. According to \cite{JCB03}, we have the following result:

\begin{prop}
The  eigenvalues of the fractional Brownian motion with Hurst exponent $H \in (0,1)$ satisfy the behavior
\begin{eqnarray*}
\lambda_p  & = &  \frac{\nu_H}{p^{2H+1}} + \mathnormal{o}\left( p^{-\frac{(2H+2)(4H+3)}{4H+5}+\delta} \right), \qquad p\gg 1 \\
\end{eqnarray*}
where $\delta > 0$ is arbitrary, $\nu_H= \frac{\mathrm{sin}(\pi H)\Gamma(2H+1)}{ \pi^{2H+1}}$, and  $\Gamma$ is the  Euler Gamma function.
\end{prop}
Therefore, when $\tau  \ll  1$, we have:
\begin{displaymath}
\lambda_p < \tau\quad \mathrm{if} \quad p > \left( \frac{ \nu_H}{\tau} \right)^{\frac{1}{2H+1}}
\end{displaymath}
We hence have the following approximation for $B_\tau$:
\begin{displaymath}
B_\tau \approx \sum_{p >  \left( \frac{ \nu_H}{\tau} \right)^{\frac{1}{2H+1}}} \frac{\nu_H}{p^{2H+1}}  + \tau  \left( \frac{ \nu_H}{\tau} \right)^{\frac{1}{2H+1}}
\end{displaymath}
Furthermore, we have:
\begin{displaymath}
 \sum_{p >  \left( \frac{ \nu_H}{\tau} \right)^{\frac{1}{2H+1}}} \frac{\nu_H}{p^{2H+1}} \approx \int_{\left( \frac{ \nu_H}{\tau} \right)^{\frac{1}{2H+1}}}^{+\infty} \frac{\nu_H}{x^{2H+1}} \, dx = \frac{\nu_H}{2H\left(\frac{ \nu_H}{\tau} \right)^{1-\frac{1}{2H+1}}}
\end{displaymath}
from which:
\begin{equation}\label{T_fbm}
B_\tau \approx C_H \tau^{1-\frac{1}{2H+1}},  \qquad \tau \ll  1
\end{equation}
where $C_H$ is a constant independent of $\tau$.

The rate of convergence for a fractional Brownian motion with Hurst parameter $H$ is $\tau^{1-\frac{1}{2H+1}}$. We note that the case $H = 1/2$ corresponds to the classical Brownian motion. We observe that the larger the Hurst parameter is (i.e. the more regular the Gaussian process is), the faster  the convergence  is. Furthermore, for $H \rightarrow 1$ the convergence rate gets close to $\tau^{2/3}$. Therefore, even for the most regular fractional Brownian motion, we are still far from the classical Monte-Carlo convergence rate.

\paragraph{Example 4 (The 1-D Mat\`ern covariance kernel)}
In this example we deal with the Mat\`ern kernel with regularity parameter $\nu > 0$ in  dimension 1:
\begin{equation}\label{mat1D}
k_{1D}(x,x';\nu,l)=\frac{2^{1- \nu}}{\Gamma(\nu)}\left( \frac{\sqrt{2\nu}|x-x'|}{l}\right)^\nu K_\nu \left( \frac{\sqrt{2\nu} |x-x'|}{l} \right)
\end{equation}
where $K_\nu$ is the modified Bessel function \cite{Ab65}. The  eigenvalues of this kernel satisfy the following asymptotic behavior \cite{Na04}:
\begin{displaymath}
\lambda_p \approx \frac{1}{p^{2 \nu}}, \qquad p \gg 1
\end{displaymath}
Following the guideline of the  Example 3 we deduce the following asymptotic behavior for $B_\tau$:
\begin{equation}\label{T_mat1D}
B_\tau \approx C_\nu \tau^{1-\frac{1}{2\nu}}, \qquad \tau \ll 1
\end{equation}
where $C_\nu$ is a constant independent of $\tau$.

This result is in agreement with the one of \cite{Ri96} who proved that for 1-dimensional kernels satisfying the Sacks-Ylvisaker of order $r$ conditions (where $r$ is an integer), the generalization error for the best linear estimator and experimental design set strategy decays as $\tau^{1-\frac{1}{2r+2}}$. Indeed, for such kernels, the eigenvalues satisfy the large-$p$ behavior  $\lambda_p \propto 1/p^{2r+2}$   \cite{R06} and by following the guideline of the previous examples we find the same convergence rate. Furthermore, our result generalizes the one of \cite{Ri96} since it provides convergence rates for more general kernels  and for any dimension (see below). Finally, our result shows that the random sampling gives the same decay rate as the optimal experimental design.

\paragraph{Example 5 (The d-D tensorised Mat\`ern covariance kernel)}
We focus here on the $d$-dimensional tensorised  Mat\`ern kernel with isotropic regularity parameter $\nu > \frac{1}{2}$.
According to \cite{PU11} the  eigenvalues of this kernel satisfy the  asymptotics:
\begin{displaymath}
\lambda_p \approx \phi(p), \qquad p\gg 1
\end{displaymath}
where the function $\phi$ is defined by:
\begin{displaymath} 
\phi(p) = \frac{\mathrm{log}(1+p)^{2(d-1)\nu}}{p^{2 \nu}} 
\end{displaymath}
Its inverse $\phi^{-1}$ satisfies:
\begin{displaymath}
\phi^{-1}(\varepsilon) = \varepsilon^{-\frac{1}{2\nu}} \left( \mathrm{log}\left(\varepsilon^{-\frac{1}{2\nu}}\right) \right)^{d-1}(1+\mathnormal{o}(1)),  \qquad  \varepsilon \ll 1
\end{displaymath}
 We hence have the  approximation:
\begin{displaymath}
B_\tau \approx    \frac{2 \nu-1}{{\phi^{-1}\left(\tau\right)}^{2\nu-1}} \mathrm{log}\left(1+\phi^{-1}\left(\tau\right)\right)^{2(d-1)\nu} + \tau\phi^{-1}\left(\tau\right)
\end{displaymath}
We can deduce the following rate of convergence for $B_\tau$:
\begin{equation}\label{T_matdD}
B_\tau \approx C_{\nu,d} \tau^{1-\frac{1}{2\nu}} \mathrm{log}\left( 1/\tau \right)^{d-1}, \qquad \tau \ll  1
\end{equation}
with $C_{\nu,d}$ a constant independent of $\tau$.

\paragraph{Example 6 (The d-D Gaussian covariance kernel)}
According to \cite{To06} the asymptotic behavior of the eigenvalues for a Gaussian kernel is:
\begin{displaymath}
\lambda_p \lesssim \mathrm{exp}\left(-p^{\frac{1}{d}}\right)
\end{displaymath}

Applying the procedure presented in the previous examples, it can be shown than the rate of  convergence of the IMSE is bounded by:
\begin{equation}\label{T_gauss}
C_d \tau \mathrm{log}\left(1/\tau\right)^d,  \qquad \tau \ll  1
\end{equation}
with $C_d$ a constant  independent of $\tau$.

\paragraph{Remark 3} We can see from  the previous examples that for smooth kernels, the convergence rate is close to $\tau$, i.e. the classical Monte-Carlo rate.

We compare the previous theoretical results on the rate of convergence of the generalization error with full numerical simulations. 
In order to observe the asymptotic convergence, we fix $n = 200$ and we consider $1/\tau$ varying from $5$ to $100$. The experimental design sets are sampled from a  uniform measure on $[0,1]$ and the observation noise is $n\tau$. To estimate the IMSE (\ref{IMSE}) we use a  trapezoidal numerical integration with 4000 quadrature points over $[0,1]$.

\begin{figure}[H]
\begin{center}
\vskip-04ex
\includegraphics[width = 6cm]{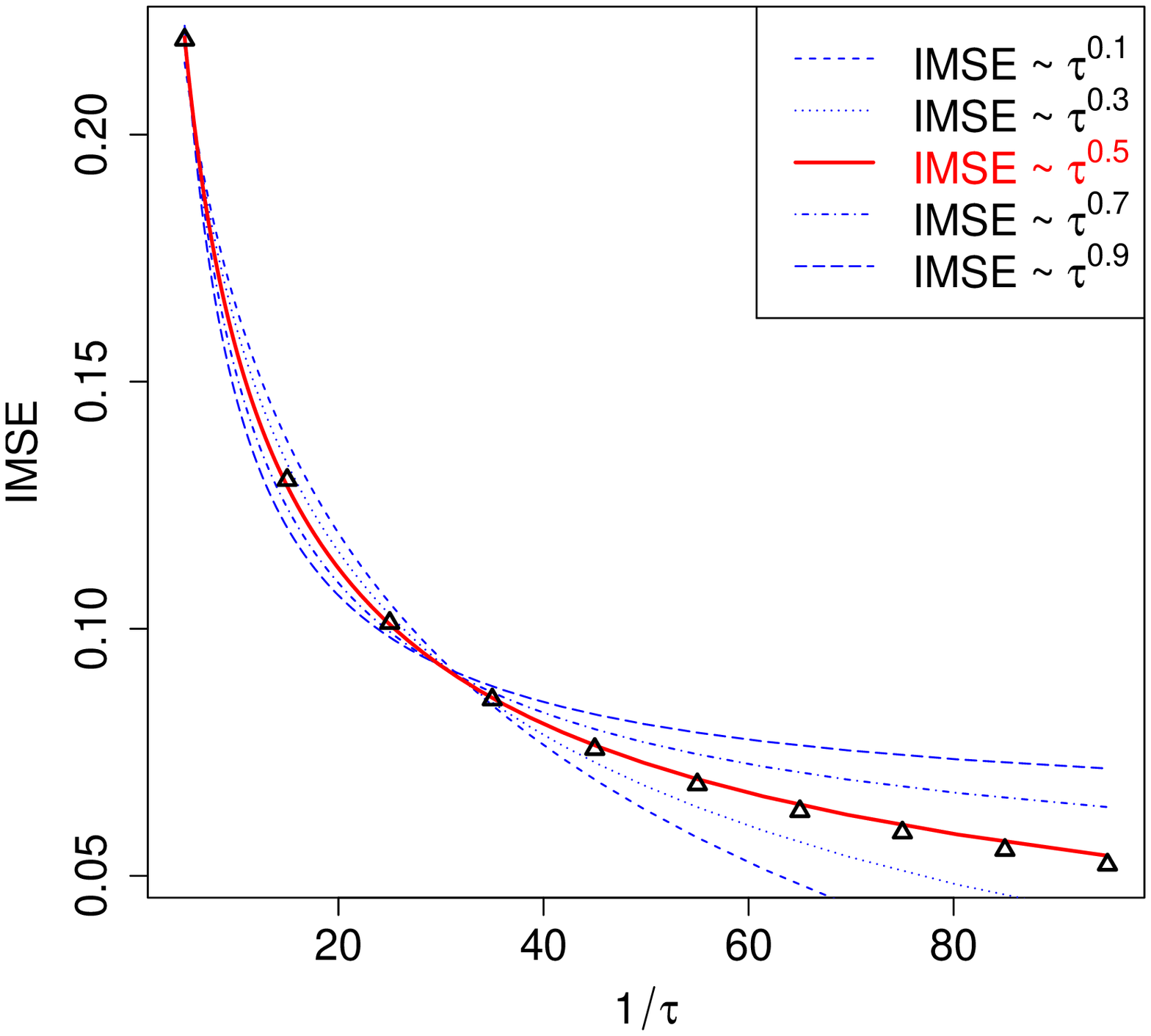}
\includegraphics[width = 6cm]{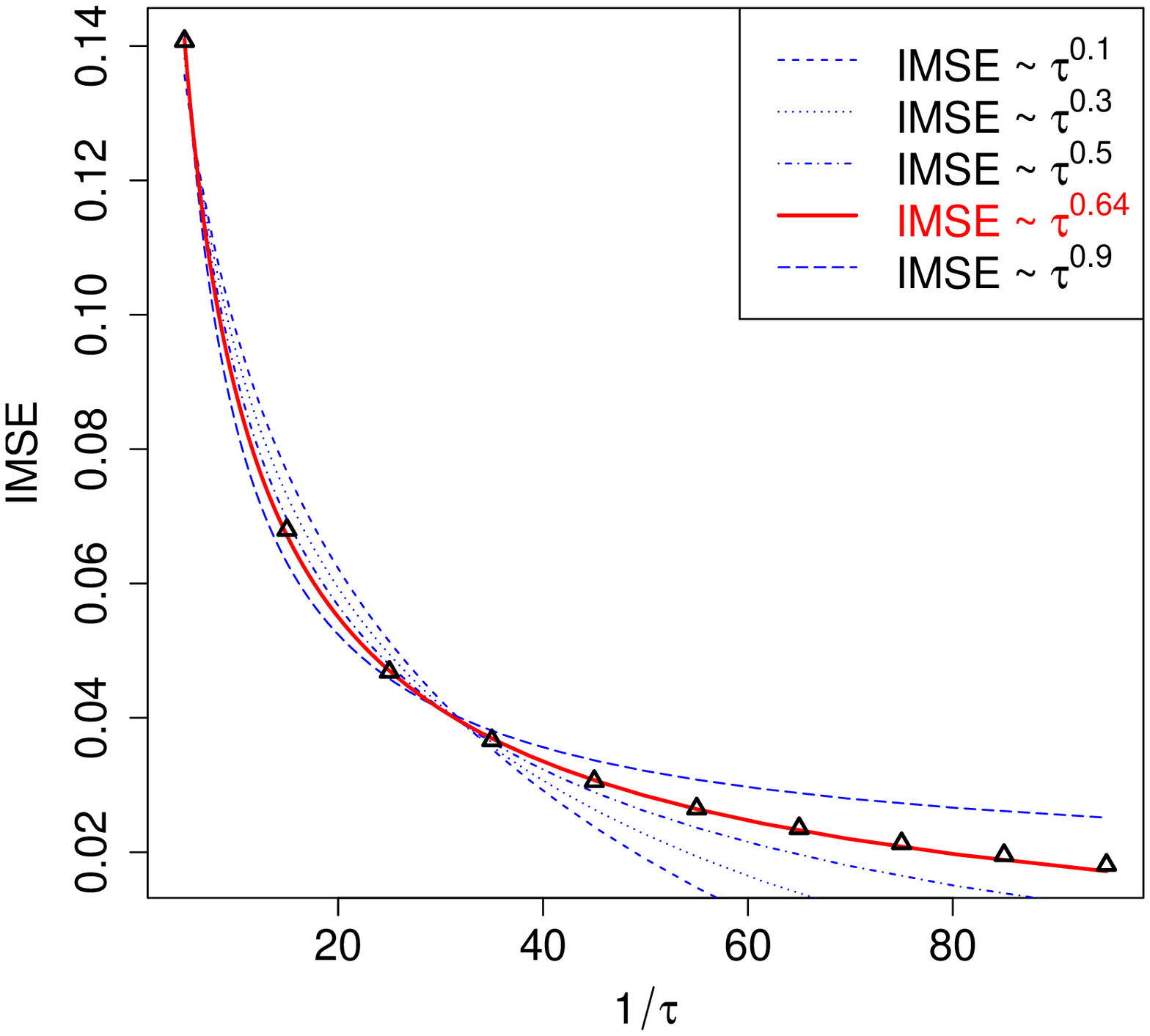}
\vskip-03ex
\caption{Rate of convergence of the IMSE  when the level of observation noise decreases for a fractional Brownian motion with Hurst parameter $H=0.5$ (left) and $H=0.9$ (right). The   number of observations is $n=200$ and the observation noise variance is $n \tau $  with $1/\tau$ varying from $5$ to $100$. The triangles represent the numerically estimated IMSE,  the solid line represents the theoretical convergence, and the other  non-solid lines represent various convergence rates.}
\label{fBm05}
\end{center}
\end{figure}

First, we deal with the 1-D fractional Brownian kernel (\ref{fBm_kernel})  with Hurst parameter $H$. We have proved that for large $n$, the IMSE decays as $\tau^{1-\frac{1}{2H+1}}$.  Figure \ref{fBm05} compares the numerically estimated convergences to the theoretical ones.

We see in Figure \ref{fBm05} that the observed rate of convergence is perfectly fitted by the theoretical one. We note that we are far from the classical Monte-Carlo rate since we are not in a non-degenerate case.  

\begin{figure}[H]
\begin{center}
\vskip-04ex
\includegraphics[width = 6cm]{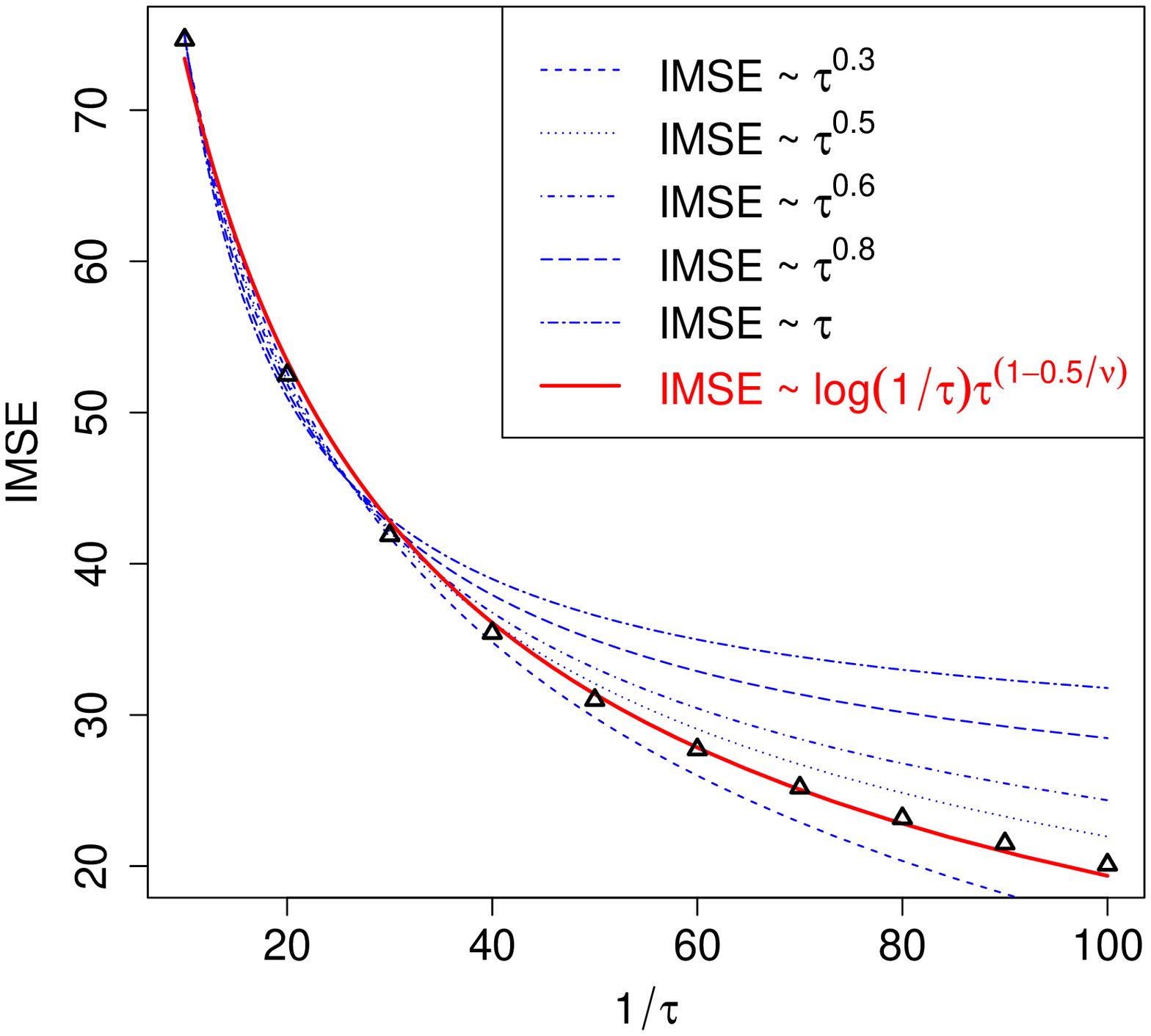}
\includegraphics[width = 6cm]{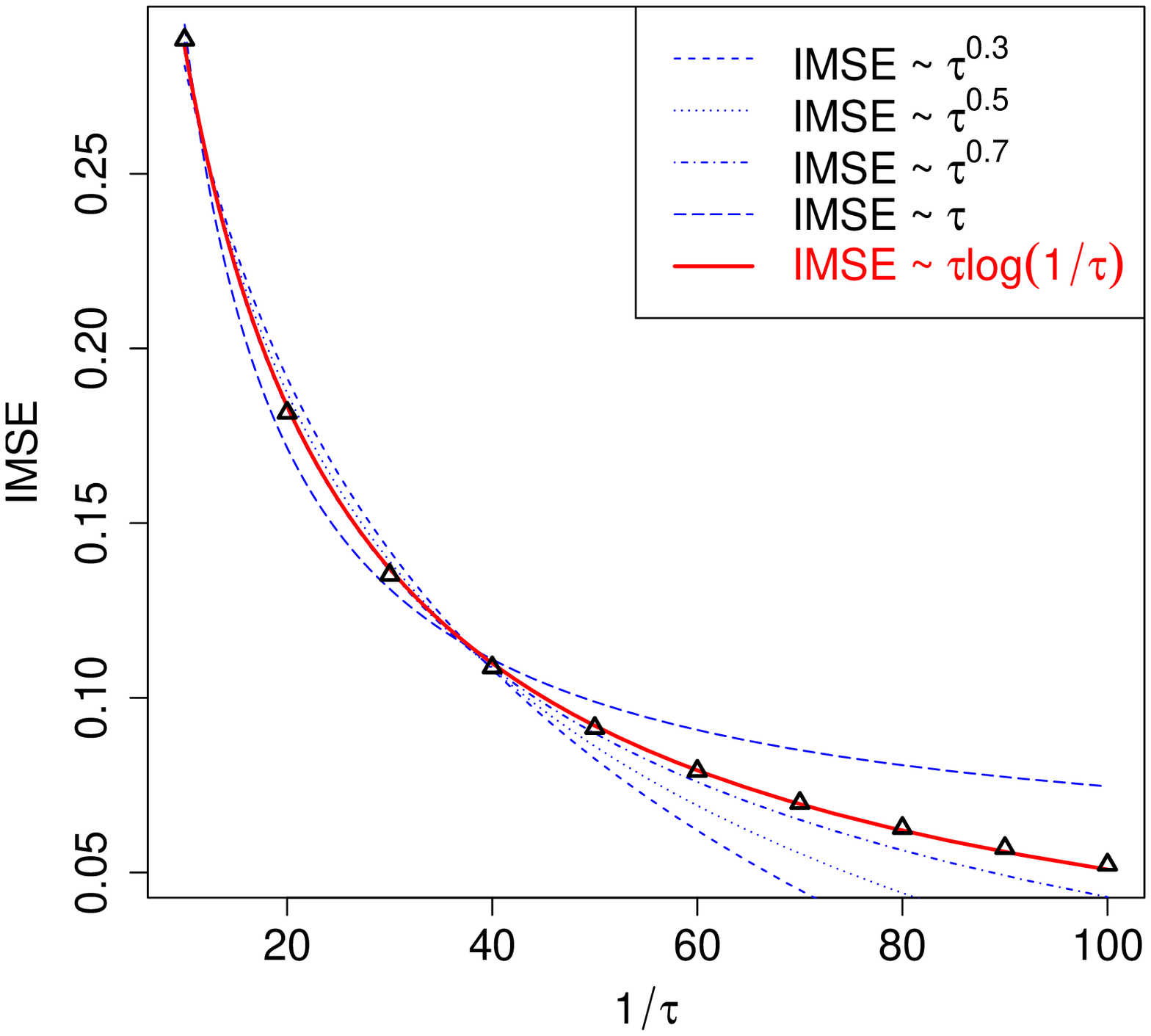}
\vskip-03ex
\caption{Rate of convergence of the IMSE  when the level of observation noise decreases for a 2-D tensorised Mat\`ern-$\frac{5}{2}$ kernel on the left hand side and for a 1-D Gaussian kernel on the right hand side. The   number of observations is $n=200$ and the observation noise variance is $n \tau $ with $1/\tau$ varying from $10$ to $100$. The triangles represent the numerically estimated IMSE,  the solid line represents the theoretical convergence, and the other non-solid lines represent various convergences.}
\label{mat2D}
\end{center}
\end{figure}

Finally, we deal with the 2-D tensorised Mat\`ern-$\frac{5}{2}$ kernel and the 1-D Gaussian kernel. The 1-dimensional Mat\`ern-$\nu$ class of covariance functions $k_{1D}(t,t';\nu,\theta) $ is given by  (\ref{mat1D}) and the 2-D tensorised Mat\`ern-$\nu$ covariance function  is given by:
\begin{equation}\label{mat2Dk}
k(x,x';\nu,\theta) = k_{1D}(x_1,x_1';\nu,\theta_1) k_{1D}(x_2,x_2';\nu,\theta_2)
\end{equation}
Furthermore, the 1-D Gaussian kernel is defined by:
\begin{displaymath}
k(x,x';\theta) = \mathrm{exp}\left(-\frac{1}{2} \frac{(x-x')^2}{\theta^2} \right)
\end{displaymath}
Figure \ref{mat2D} compares the numerically observed convergence of the IMSE to the theoretical one when $\theta_1 = \theta_2 = 0.2$ for the Mat\`ern-$\frac{5}{2}$ kernel and when $\theta = 0.2$ for the Gaussian kernel. We see in figure \ref{mat2D} that the theoretical rate of convergence is a sharp approximation of the observed one. 

\section{Applications of the learning curve}\label{allocation}

Let us consider that we want to approximate the function $x\in \mathbb{R}^d \rightarrow f(x) $ from  noisy observations at fixed  points $(x_i)_{i=1,\dots,n}$, with $n \gg 1$, sampled from the design measure $\mu$ and with $s_i$ replications at each point   $x_i$.

In this section, we consider the situation described in Example \ref{ex1}:
\begin{itemize}
\item The budget $T$ is defined as the sum of repetitions on all points of the experimental design set - i.e.  $T = \sum_{i=1}^{n} s_i$.
\item An observation $z_i^n$ at point $x_i$   has a noise variance equal to $\sigma_\varepsilon^2(x_i)/s_i$ with $i=1,\dots,n$.
\end{itemize}

In Subsection \ref{EsimationT} we present how to determine the needed budget $T$ to achieve a prescribed precision. Then,  in Subsection \ref{AllocationT},  we address the problem of the optimal allocation   $\{s_1, s_2, \dots, s_{n}\}$ for a given budget $T$.

\subsection{Estimation of the budget required to reach a prescribed precision}\label{EsimationT}

Let us consider a prescribed  generalization error  denoted by  $\boldsymbol{\bar{\varepsilon}}$. The purpose of this subsection is to determine from an initial budget $T_0$ the budget $T$ for  which the generalization error reaches   the value  $\boldsymbol{\bar{\varepsilon}}$.
We handle this issue by considering a uniform allocation $s_i = s$   with  $i=1,\dots,n$ and  a constant reduced noise variance $\sigma_\varepsilon^2$.

First, we build an initial experimental design set $(x_i^\mathrm{train})_{i=1,\dots,n}$ sampled with respect to the design measure $\mu$ and with $s^*$ replications at each point  such that $T_0=ns^*$. From the $s^*$ replications $(z_{i,j})_{j=1,\dots,s^*}$,  we can estimate the observation noise variances   $\sigma_\varepsilon^2(x_i^\mathrm{train})$ with a classical empirical estimator: $ \sum_{j=1}^{s^*}(z_{i,j}-z_i^n)^2/(s^*-1) ,\, z_i^n = \sum_{j=1}^{s^*}z_{i,j}/s^*$. Then, we consider a constant reduced noise variance $\sigma_\varepsilon^2$ equal to the  mean $\int_{\mathbb{R}^d} \sigma_\varepsilon^2(x )\,d\mu(x)$ estimated with  $\sum_{i=1}^n \sigma_\varepsilon^2(x_i^\mathrm{train}) /n$. 

Second, we use the observations $z_i^n = (\sum_{j=1}^{s^*}z_{i,j})/s^*$ to estimate the covariance kernel $k(x,x')$. In practice, we consider a parametrized family of covariance kernels and we select the parameters which maximize the likelihood \cite{S99}.

Third, from  Proposition \ref{theo3} we can get the expression of the generalization error decay with respect to $T$  (denoted by  IMSE$_T$). Therefore, we just have to determine the budget $T$ such that $\mathrm{IMSE}_T = \boldsymbol{\bar{\varepsilon}}$. In practice, we will not use Proposition \ref{theo3} but the asymptotic results described in Section \ref{rate}.

This strategy will be applied to an industrial case in Section \ref{application}. We note that in the application presented  in Section \ref{application}, we have $s^*=1$. In fact, in this example  the observations are themselves obtained by an empirical mean of a Monte-Carlo sample and thus the noise variance can be estimated without processing replications.

\subsection{Optimal resource allocation for a given budget}\label{AllocationT}

In this subsection, we consider a  fixed budget $T$. As presented in Subsection \ref{EsimationT}, to determine this budget  we make the approximation of a reduced noise variance $\sigma_\varepsilon^2(x)$ independent of $x$ and we consider the uniform allocation $s_i=s$.

Despite the fact that   the uniform allocation   $s_i = s$  are needed to determine  $T$, in order to provide the optimal resource allocation - i.e. the sequence of integers $\{s_1, s_2, \dots, s_{n}\}$  minimizing the generalization error - it is worth taking into account the heterogeneity of the noise. For a Monte-Carlo based simulator, the number of repetitions $s$ could represent the number of MC particles and the  procedure presented below can be applied.

Determining the optimal allocation of the budget $T$ whatever the Gaussian process for a heterogeneous noise is an open and non-trivial problem.
 To solve this problem,  we first consider the continuum approximation in which we look for an optimal sequence of real numbers  $(s_i)_{i=1,\dots,n}$  and then  we  round the optimal solution to obtain a quasi-optimal integer-valued allocation $(s_{i,\mathrm{int}})_{i=1,\dots,n}$.
The following proposition gives the optimal resource allocation under certain restricted conditions for the continuous case. The reader is referred to \cite{MZ11} for a proof of this proposition in a different framework (the proof uses the Karush-Kuhn-Tucker  approach to solve the minimization problem with equality and inequality constraints). We note that the optimal allocation given in Proposition \ref{theo5} for a fixed budget $T$ can also be used for any $n > 0$  and for any experimental design sets. In particular, it  is not restricted to the case $n$ large.

\begin{prop}\label{theo5}
Let us consider $Z(x)$ a Gaussian process with a known mean and covariance kernel $k(x,x') \in \mathcal{C}^0(\mathbb{R}^d \times \mathbb{R}^d)$ with $\sup_x k(x,x) < \infty$. Let  $(x_i)_{i=1,\dots,n}$ be a given  experimental design set of $n$ points sorted such  that the sequence $\left(\frac{k(x_j,x_j)+\sigma_\varepsilon^2(x_j)}{\sqrt{c(x_j) \sigma_\varepsilon^2(x_j)}}\right)_{j=1,\dots,n}$ is non-increasing, where $\sigma_\varepsilon^2(x_i)$ is the reduced  noise variance of an observation  at point $x_i$,  $c(x) = \int_{\mathbb{R}^d}k(x',x)^2  \, d\eta(x') $ and $\eta(x)$ is  a positive measure used to calculate the IMSE.  When the covariance matrix $K$ is diagonal, the real-valued allocation $(s_i)_{i=1,\dots,n}$ minimizing the generalization error:
\begin{equation}\label{IMSE_DELTA}
\mathrm{IMSE} = \int_{\mathbb{R}^d}{\left(k(x,x) - k(x)^T (K+ \Delta)^{-1}k(x) \right) \, d\eta(x)}
\end{equation}
under the constraints $\sum_{i=1}^{n} s_i = T$ and $s_i \geq 1, \, \forall i =1,\dots,n$  is given by:
\begin{equation}
s_i^{\mathrm{opt}} = \left\{
\begin{array}{ll}
1  & i \leq i^*\\
\frac{1}{k(x_i,x_i)}\left(
 \frac{ \sqrt{c(x_i)\sigma_\varepsilon^2(x_i)}}{\sum_{j  = i^*+1}^{n}\frac{\sqrt{c(x_j)\sigma_\varepsilon^2(x_j)}}{k(x_j,x_j)}}\left( T - i^*+   \sum_{j = i^*+1}^{n}\frac{\sigma_\varepsilon^2(x_j)}{k(x_j,x_j)}\right)
 -  \sigma_\varepsilon^2(x_i)\right) &
i  > i^*\\
\end{array}
\right.
\end{equation}
where   $ \Delta = \mathrm{diag}\left[\left( \frac{\sigma_\varepsilon^2(x_i)}{s_i} \right)_{i=1,\dots,n}\right]$and:
\begin{equation}
i^* = \max \left\{ i=1,\dots,n \quad \mathrm{such \, that} \quad \frac{k(x_i,x_i)+\sigma_\varepsilon^2(x_i)}{\sqrt{c(x_i) \sigma_\varepsilon^2(x_i)}}
\geq
\frac{T-i+\sum_{j=i+1}^{n} \frac{\sigma_\varepsilon^2(x_j)}{k(x_j,x_j)}}{\sum_{j=i+1}^{n}\frac{\sqrt{c(x_j) \sigma_\varepsilon^2(x_j)}}{k(x_j,x_j)}}
\right\}
\end{equation}
By convention, if:
\begin{equation}
  \frac{k(x_i,x_i)+\sigma_\varepsilon^2(x_i)}{\sqrt{c(x_i)\sigma_\varepsilon^2(x_i)}}
<
\frac{T-i+\sum_{j=i+1}^{n} \frac{\sigma_\varepsilon^2(x_j)}{k(x_j,x_j)}}{\sum_{j=i+1}^{n}\frac{\sqrt{c(x_j) \sigma_\varepsilon^2(x_j)}}{k(x_j,x_j)}},\qquad \forall i=1, \dots,n
\end{equation}
then $i^*=0$.
\end{prop}

The optimization problem in  Proposition \ref{theo5} admits a solution if and only if $T \geq n$ which  reflects the fact that  $n$  simulations are already available. Furthermore, when $T$ is large enough, we have $i^*=0$ and the solution has the following form:
\begin{equation}\label{riopt}
s_i^{\mathrm{opt}}  = \frac{1}{k(x_i,x_i)}\left(
 \frac{ \sqrt{c(x_i)\sigma_\varepsilon^2(x_i)}}{\sum_{j  = 1}^{n}\frac{\sqrt{c(x_j)\sigma_\varepsilon^2(x_j)}}{k(x_j,x_j)}}\left( T +  \sum_{j = 1}^{n}\frac{\sigma_\varepsilon^2(x_j)}{k(x_j,x_j)}\right)
 -  \sigma_\varepsilon^2(x_i)\right) 
\end{equation}

While  Proposition \ref{theo5} gives a continuous optimal allocation, an admissible allocation must be an integer-valued sequence. Therefore, as mentioned previously, we solve the optimization problem with the continuous approximation and then we round the continuous solution to obtain a  quasi-optimal integer-valued  solution $s_{i,\mathrm{int}}^{\mathrm{opt}}$. 
The rounding is performed by solving the following problem:

Find $J$ such that $\sum_{i=1}^{n} s_{i,\mathrm{int}}^{\mathrm{opt}} = T$ with:
\begin{displaymath}
s_{i,\mathrm{int}}^{\mathrm{opt}} = \left\{
\begin{array}{ll}
\big[ s_i^{\mathrm{opt}} \big] +1 &  i \leq J \\
\big[  s_i^{\mathrm{opt}} \big] & i > J
\end{array}
\right.
\end{displaymath}
where $[ x ]$ denotes  the integer part of a real number $x$.

We note that this allocation is not optimal in general (i.e. when $K$ is not diagonal). Nevertheless we have numerically observed that it remains efficient in  general cases and is   better than the uniform allocation strategy. We note that the numerical comparison has been performed with different kernels (Gaussian, Mat\`ern-$\frac{5}{2}$, Mat\`ern-$\frac{3}{2}$, exponential, Brownian and triangular \cite{R06})  and in dimension one and two with a number of observations varying between 10 and 400. Furthermore, two  types of experimental design sets have been tested, one is a random set sampling from the uniform distribution and the other one is a regular grid.

 Proposition \ref{theo5} shows that it is worth  allocating  more resources at locations where the reduced noise variance $\sigma_\varepsilon^2(x)$ and  the   quantity $c(x_i) = \int_{\mathbb{R}^d}k(x,x_i)^2  \, d\eta(x) $  (representing the local concentration of the IMSE)  are more  important.  

\section{Industrial Case: code MORET}\label{application}

We illustrate in this section an industrial application of our results about the rate of convergence of the IMSE. The case  is about  the  safety assessment of a nuclear system containing fissile materials. The system is modeled by a neutron transport code called MORET \cite{Fe05}.  In particular, we study a benchmark system of dry $PuO_2$ storage.  We note that we are in the framework presented in Example \ref{ex1}.

This section is divided into 3 parts. First, we present the Gaussian process regression model built on an initial experimental design set. Then we apply the strategy described in Section \ref{EsimationT} to determine the computational budget $T$ needed to achieve a prescribed precision.  Finally, we allocate the resource $T$ on the experimental design set.

\subsection{Data presentation}

The benchmark system safety  is evaluated through the neutron multiplication factor $k_{\mathrm{eff}}$. This is  our output of interest that we want to surrogate. This factor models the criticality of a chain nuclear reaction:
\begin{itemize}
\item $k_{\mathrm{eff}} > 1$ leads to an uncontrolled chain reaction due to an increasing neutron population.
\item $k_{\mathrm{eff}} = 1$ leads to a self-sustained chain reaction with a stable neutron population.
\item $k_{\mathrm{eff}} < 1$ leads to a faded chain reaction due to an decreasing neutron population.
\end{itemize}
The neutron multiplication factor depends on many parameters and it is evaluated using the stochastic simulator called MORET.
 We focus here on two parameters:
\begin{itemize}
\item $d_{\mathrm{PuO}_2} \in [0.5 , 4] \mathrm{g.cm}^{-3}$, the density of the fissile powder. It is scaled in this section  to $[0,1]$.
\item $d_{\mathrm{water}} \in [0,1] \mathrm{g.cm}^{-3}$, the density of water between storage tubes.
\end{itemize}
The other parameters are fixed to a nominal value given by an expert and we use the notation $x =  (d_{\mathrm{PuO}_2}, d_{\mathrm{water}})$ for the input parameters.\\

The MORET code provides outputs of the following form:
\begin{displaymath}
k_{\mathrm{eff},s}(x) = \frac{1}{s} \sum_{j=1}^{s}Y_j(x)
\end{displaymath}
where $(Y_j(x))_{j=1,\dots,s}$ are  realizations of  independent and identically distributed random variables  which are themselves obtained by an empirical mean of  a Monte-Carlo sample of 4000 particles. From these particles, we can also estimate the variance $\sigma^2_{\varepsilon}(x)$ of the observation $Y_j(x)$ by a classical empirical estimator. The simulator gives noisy  observations and the variance of an observation $k_{\mathrm{eff},s}(x)$ equals $\sigma^2_\varepsilon(x)/s$. 

A  large data base $(Y_j(x_i))_{ i=1,\dots,5625, j=1,\dots,200}$ is available to us. We divide  it into a training set and a test set. Let us denote by $Y_j(x_i)$ the j$^\mathrm{th}$ observation at point $x_i$ - the 5625 points $x_i$ of the data base come from a $75 \times 75$ grid over $[0,1]^2$. The training set consists of $n= 100$ points $(x_i^{\mathrm{train}})_{i=1,\dots,n}$ extracted from the complete data base using a maximin LHS  and of the first observations $(Y_1(x_i^{\mathrm{train}}))_{i=1,\dots,100}$. We will use the other 5525  points as a test set.

The aim of the study is - given the training set -  to predict the budget needed to achieve a prescribed precision for the surrogate model and to allocate optimally these resources. More precisely, let us denote by $s_i$ the resource allocated to the point $x_i^\mathrm{train}$ of the experimental design set. First, we want to determine the budget $T = \sum_{i = 1}^{n} s_i$ which allows us to achieve the target precision (see Subsection \ref{EsimationT}). Second, we want to determine the best resource allocation  $(s_i)_{i=1,\dots, n}$ (see Subsection \ref{AllocationT}).

To evaluate the needed computational budget $T$ the observation noise variance  $\sigma^2_{\varepsilon}(x)$ is approximated by  a constant $\bar{\sigma}^2_\varepsilon$ in order to fit with the hypotheses of the theorem. The constant variance equals the mean $ \int_{\mathbb{R}^2}\sigma^2_{\varepsilon}(x)\,d\mu(x) $ of the noise variance   which is here estimated by  $  \bar{\sigma}^2_\varepsilon = \frac{1}{100} \sum_{i=1}^{100}\sigma^2_{\varepsilon}(x_i^{\mathrm{train}})  = 3.3.10^{-3}$. Furthermore, we look for a uniform budget allocation, i.e. $s_i = s$ $\forall i=1,\dots,n$. In this case, the total computational budget is $T = ns$.

\subsection{Model selection}

To build the model, we consider the training set plotted in figure \ref{Nopt1}. It is composed of the $n=100$ points $(x_i^{\mathrm{train}})_{i=1,\dots,n}$ which are uniformly spread on $Q = [0,1]^2$.

Let us suppose that the response is a realization of a Gaussian process with a tensorised Mat\`ern-$\nu$ covariance function. The 2-D  tensorised Mat\`ern-$\nu$ covariance function $k(x,x';\nu,\theta)$  is given in (\ref{mat2Dk}).
 The hyper-parameters are estimated by maximizing the concentrated Maximum Likelihood \cite{S99}:
\begin{displaymath}
-\frac{1}{2}(z-m)^T(\sigma^2K+\sigma^2_{\varepsilon}I)^{-1}(z-m)-\frac{1}{2} \mathrm{det}(\sigma^2K+\bar{\sigma}^2_\varepsilon I)
\end{displaymath}
where $K = [k(x_i^\mathrm{train},x_j^\mathrm{train};\nu,\theta) ]_{i,j=1,\dots,n}$, $I$ is the identity matrix,  $\sigma^2$ the variance parameter, $m$ the mean of $k_{\mathrm{eff},s}(x) $  and $z = (Y_1(x_1^\mathrm{train}), \dots, Y_1(x_n^\mathrm{train}))$  the observations at points in the training set. The mean of $k_{\mathrm{eff},s}(x) $ is estimated by $m=\frac{1}{100} \sum_{i=1}^{100} Y_{1}(x_i^{\mathrm{train}}) = 0.65$. 

Due to the fact that the convergence rate is strongly dependent of the regularity parameter $\nu$, we have to perform a good estimation of this hyper-parameter to  evaluate the model error  decay accurately. Note that we cannot have a closed form expression for the estimator of $\sigma^2$, it hence has to be  estimated jointly  with $\theta$ and $\nu$.

Let us consider the vector of parameters $\phi = (\nu, \theta_1, \theta_2, \sigma^2)$. In order to perform the maximization, we have first randomly generated a set of 10,000 parameters $(\phi_k)_{k=1,\dots,10^4}$  on the domain $[0.5,3] \times[0.01,2] \times [0.01,2]\times [0.01,1]$. We have then selected the 150 best parameters (i.e. the ones maximizing the concentrated Maximum Likelihood) and we have started a quasi-Newton based maximization from these parameters. More specifically, we have used the BFGS method \cite{Sh70}. Finally, from the results of the 150 maximization procedures, we have selected the best parameter. We note that the quasi-Newton based maximizations have all converged to two parameter values, around 30\% to the actual maximum and 70\% to another local maximum.   

The estimation of the hyper-parameters are $\nu = 1.31$,  $\theta_1 = 0.67$, $\theta_2 = 0.45$ and $\sigma^2= 0.24 $. This means that we have a rough surrogate model which is not differentiable and  $\alpha$-H\"older continuous with exponent $\alpha = 0.81$.
The variance of the observations  is   $\bar{\sigma}^2_\varepsilon=3.3.10^{-3}$, using the same notations as  Example \ref{ex1}, we have $\tau = \bar{\sigma}^2_\varepsilon/T_0  $ with $T_0 = n$ (it corresponds to $s=1$).

The IMSE of the Gaussian process regression  is $\mathrm{IMSE}_{T_0}=1.0.10^{-3}$ and its empirical mean squared error is $\mathrm{EMSE}_{T_0} = 1.2.10^{-3}$ . To compute the empirical mean squared error (EMSE), we use the observations $(Y_j(x_i))_{i=1,\dots,5525,\,  j=1   \dots, 200}$ with $x_i \neq x_k^{\mathrm{train}}$ $\forall k=1,\dots,100, i=1,\dots, 5525$ and to compute the IMSE (\ref{IMSE}) (that depends only on the positions of the training set and on the selected hyper-parameters) we use a trapezoidal numerical integration into a $75 \times 75$ grid over $[0,1]^2$. For $s=200$, the observation  variance of the output $k_{\mathrm{eff},s}(x)$   equals $\frac{ \bar{\sigma}^2_\varepsilon}{200}=1.64.10^{-5}$ and  is neglected for the estimation of the empirical error. We can  see that the IMSE is  close to the empirical mean squared error which means that our model  describes the observations accurately. 

\subsection{Convergence of the  IMSE}

According to (\ref{T_matdD}), we have the following convergence rate for the IMSE:

\begin{equation}
\mathrm{IMSE} \sim \mathrm{log}(1/\tau) \tau^{1-\frac{1}{2\nu}} = \frac{\mathrm{log}(T/\bar{\sigma}^2_{\varepsilon})}{(T/\bar{\sigma}^2_{\varepsilon})^{1-\frac{1}{2\nu}}}
\end{equation}
where the model parameter $\nu$ plays a crucial role. We can therefore expect that the IMSE decays as (see Subsection \ref{EsimationT}):
\begin{equation}\label{IMSE_T}
\mathrm{IMSE}_T =  \mathrm{IMSE}_{T_0}  \frac{\mathrm{log}(T/\bar{\sigma}^2_{\varepsilon})}{(T/\bar{\sigma}^2_{\varepsilon})^{1-\frac{1}{2\nu}}}/\frac{\mathrm{log}(T_0/\bar{\sigma}^2_{\varepsilon})}{(T_0/\bar{\sigma}^2_{\varepsilon})^{1-\frac{1}{2\nu}}}
\end{equation}

\begin{figure}[H]
\begin{center}
\vskip-04ex
\includegraphics[width =7cm]{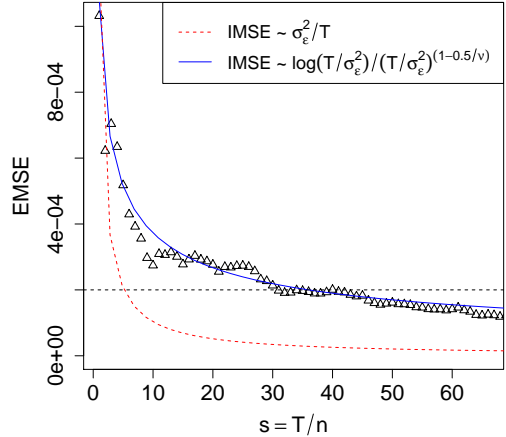}
\vskip-03ex
\caption{Comparison between Empirical mean squared error (EMSE) decay and theoretical IMSE decay for $n=100$ when the total budget $T=ns$ increases.  The triangles represent the Empirical MSE, the solid line represents the theoretical decay,  the horizontal dashed line represents the desired accuracy and the dashed  line the classical M-C convergence. We see that Monte-Carlo decay does not match the empirical MSE and it is too fast. }
\label{empdecay1}
\end{center}
\end{figure}

Let us assume that we want to reach an IMSE of $\boldsymbol{\bar{\varepsilon}} = 2.10^{-4}$. According to the IMSE decay and the fact that the IMSE for the budget $T_0$ has been estimated to be equal to $1.0.10^{-3}$, the total budget required is $T = ns = 3600$, \emph{i.e.} $s = 36$.  
Figure \ref{empdecay1} compares the empirical mean squared error convergence and the predicted convergence (\ref{IMSE_T}) of the IMSE.

We see empirically that the EMSE of $\boldsymbol{\bar{\varepsilon}} = 2.10^{-4}$ is achieved  for $s = 31$.  This shows that the predicted IMSE and the empirical MSE are very close and that the selected kernel captures  the regularity of the response accurately.

Let us consider the classical Monte-Carlo convergence rate $ \bar{\sigma}_\varepsilon^2/T$, which  corresponds to the convergence rate of degenerate kernels, \emph{i.e.} in the finite -dimensional case. Figure \ref{empdecay1} compares the theoretical rate of convergence of the IMSE with the classical Monte-Carlo one. We see that the Monte-Carlo decay is too fast and does not represent correctly the empirical MSE decay. If we had considered the rate of convergence $\mathrm{IMSE} \sim  \bar{\sigma}_\varepsilon^2/T$, we would have reached an  IMSE of $\boldsymbol{\bar{\varepsilon}}  = 2.10^{-4}$ for $s = 6$ (which is very far from the observed value $ s = 31$).

\subsection{Resources allocation}

We have determined  in the previous section the computational budget required to reach an IMSE of $2.10^{-4}$. We observe that the predicted allocation is accurate since it gives an empirical MSE close to $2.10^{-4}$. To calculate the observed MSE, we uniformly allocate the computational budget on  the points of the training set. We know that this allocation is optimal when the variance of the observation noise is homogeneous. Nevertheless, we are not in this case and to build the final model we allocate the budget taking into account the heterogeneous noise level $\sigma^2_{\varepsilon}(x)$.   
We note that the  total budget is $T = \sum_{i = 1}^{n} s_i$ where $n=100$ is the number observations  and $s_i$ the budget allocated to the point $x_i^\mathrm{train}$.

From (\ref{riopt}), when the input parameter distribution $\mu$ is uniform on $[0,1]$  and  for a diagonal covariance matrix, the optimal allocation is given by:
\begin{equation}\label{Nopt}
s_i = \frac{1}{\sigma^2}\left(
 \frac{ \sqrt{\sigma^2_\varepsilon(x_i)}}{\sum_{j=1}^{n} \sqrt{\sigma^2_\varepsilon(x_j)}}\left( \sigma^2 T +  \sum_{j=1}^{n}\sigma^2_\varepsilon(x_j)\right)
 -   \sigma^2_\varepsilon(x_i)\right)
\end{equation}
 
\begin{figure}[H]
\begin{center}
\vskip-04ex
\includegraphics[width = 6cm]{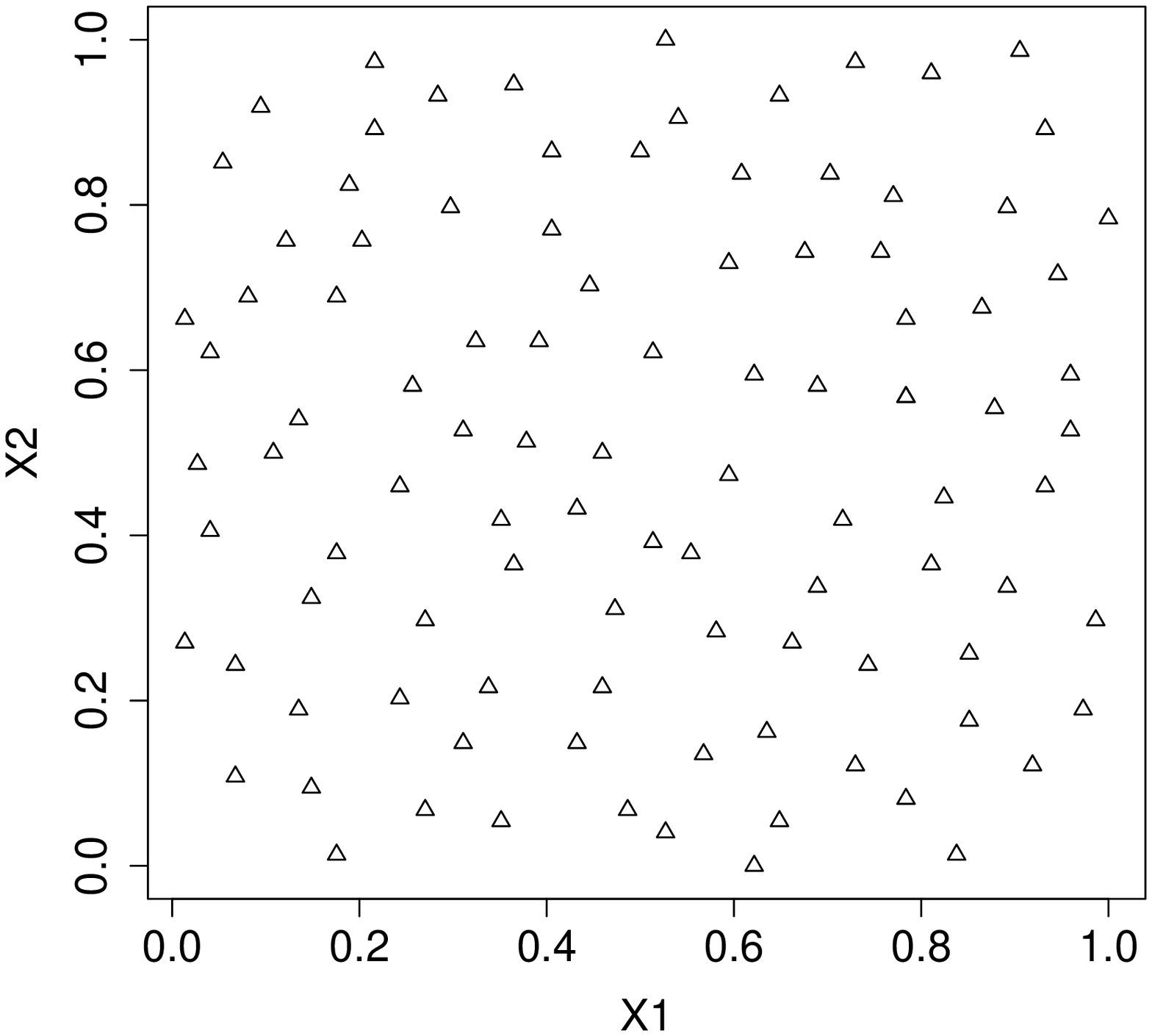}
\includegraphics[width =6cm]{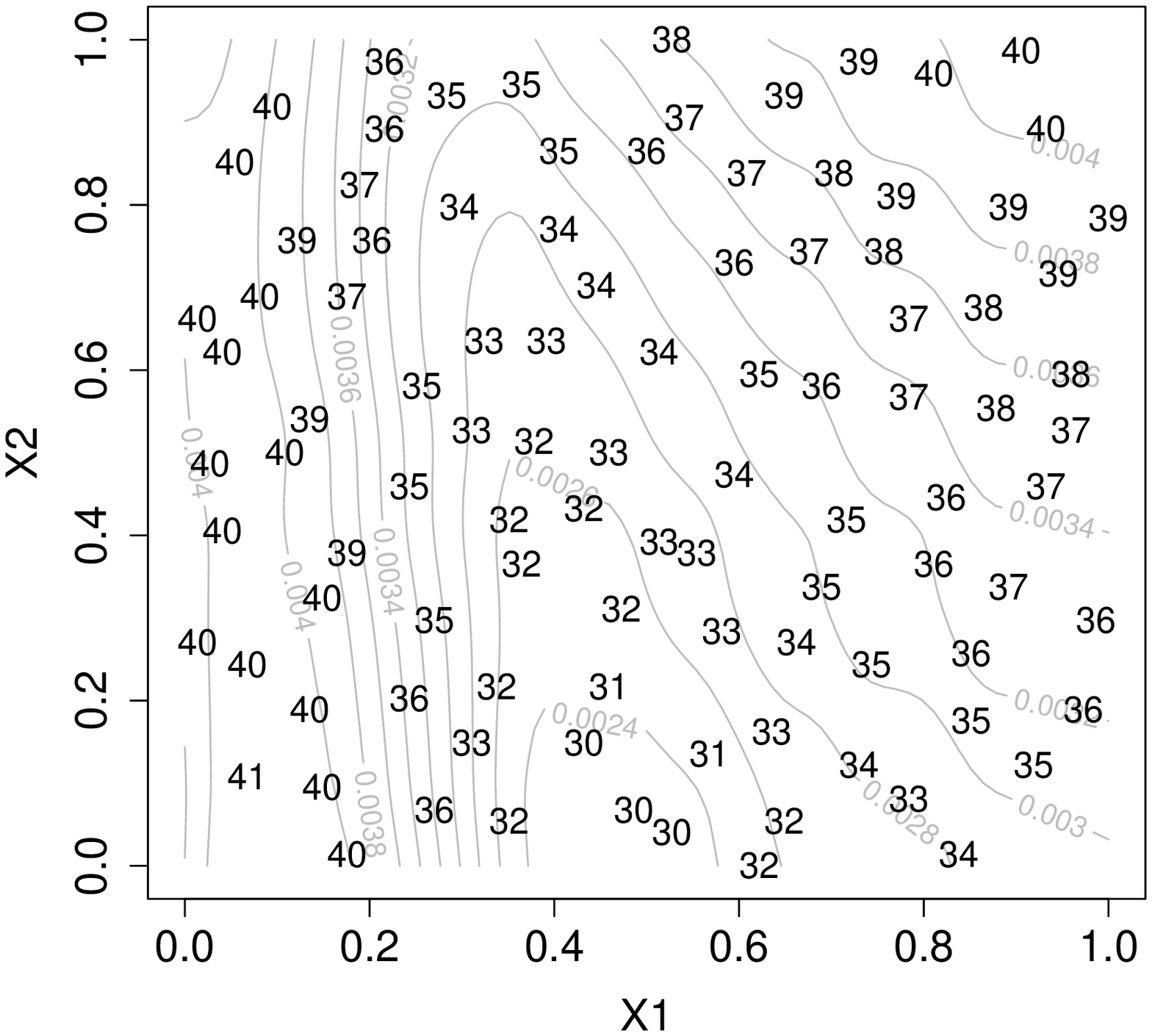}
\vskip-03ex
\caption{On the left hand side: initial experimental design set with $n=100$. On the right hand side: noise level dependence of the resources allocation. The solid lines represent the reduced noise variance $\sigma_\varepsilon^2(x)$ contour plot and the numbers represent the resources $(s_i)_{i=1,\dots,n}$ allocated to the points of the experimental design set.}
\label{Nopt1}
\end{center}
\end{figure}

 Here we  use this allocation to build the model. Let us consider that we do not have  observed the empirical MSE decay, we hence consider the budget given by the theoretical decay $T = 3600$. The allocation given by  equation (\ref{Nopt}) after the rounding procedure is illustrated in figure \ref{Nopt1} with the contour of the noise level.

We see in figure \ref{Nopt1} that the resources allocation is more important at points where the noise variance is higher. Table \ref{resalloc1} compares the performances of the two models build with the  two allocations on the test set.

\begin{table}[H]
\begin{center}
\begin{tabular}{|c|c|c|}
\hline
  & Uniform Allocation & Optimal Allocation \\
\hline
MSE & $1.94.10^{-4}$ & $1.86.10^{-4}$ \\
\hline
MaxSE & $3.66.10^{-2}$ & $3.38.10^{-2}$  \\
\hline
\end{tabular}
\end{center}
\caption{  Comparison between uniform and optimal (under the condition $K$ diagonal) allocation of resources. }
\label{resalloc1}
\end{table}

We see in Table \ref{resalloc1} that the budget allocation given by the equation (\ref{Nopt}) gives predictions slightly more accurate than the uniform one.

\section{Acknowledgments}

The authors are grateful to  Dr. Yann Richet of the IRSN - Institute for Radiological Protection and Nuclear Safety - for providing the data for the industrial case through the reDICE project.

\section{Conclusion}

The main result of this paper is a theorem giving the Gaussian process regression mean squared error  when the number of observations is large and the observation noise variance is   proportional to the number of observations. 
The asymptotic value of the mean squared error  is derived in terms of the eigenvalues and eigenfunctions of the   covariance function and holds for degenerate and non-degenerate kernels and for any dimension. We emphasize that a noise variance proportional to the number of  observations is natural in the framework of experiments with replications or Monte-Carlo simulators. 

From this theorem, we can deduce the asymptotic behavior  of  the generalization error - defined in this paper as the Integrated Mean Squared Error - as a function of the reduced observation noise variance (it corresponds to the noise variance  when the number of observations equals one). This result generalizes previous ones which give this behavior in dimension one or two or for a restricted class of covariance kernels (for degenerate ones).
 The significant differences between the rate of convergence of degenerate and non-degenerate kernels highlight the relevance of our theorem which holds for non-degenerate kernels. This is especially important as usual kernels for Gaussian process regression are non-degenerate.

Our work  deals   with Gaussian process regression when the variance of the noise can be reduced by increasing the budget (i.e. the number of replications at each point). Our results are of practical interest in this case since it gives the total budget needed to reach a precision prescribed by the user. Nonetheless, it holds under the assumptions of homoscedastic observation noise. Despite the fact that this assumption is relevant to evaluate the budget, it is not optimal to determine the resources allocation. Indeed,  in this case  it is worth  taking into account the noise variance heterogeneity and using a non-uniform allocation. We describe the resulting error reduction under restricted conditions. We have observed on test cases that our non-uniform allocation is better than the uniform one in more general cases although it is not optimal anymore.

\appendix

\section{Proof of the main theorem}\label{A_proof}

\subsection{Proof of Theorem \ref{theo1}: the degenerate case}

 The  proof in the degenerate case follows the lines of the ones  given by \cite{OV99}, \cite{R06}  and  \cite{VP09}. For a degenerate kernel, the number $\bar{p}$ of non-zero eigenvalues is finite. Let us denote $
\Lambda = \mathrm{diag}(\lambda_i)_{1 \leq i \leq \bar{p}}$,  $ \phi(x) = (\phi_1(x), \dots, \phi_{\bar{p}}(x))$ and $ \Phi = \left(\begin{array}{ccc} \phi(x_1)^T  &  \dots  &  \phi(x_{n})^T \end{array} \right)^T$.
The MSE of the Gaussian process regression  is given by:
\begin{displaymath}
\sigma^2(x) = \phi(x) \Lambda \phi(x)^T  -  \phi(x) \Lambda \Phi^T \left(\Phi \Lambda \Phi^T + n\tau I \right)^{-1} \Phi \Lambda \phi(x)^T
\end{displaymath}

Thanks to the Woodbury-Sherman-Morrison formula  and according to \cite{OV99} and  \cite{VP09} the Gaussian process regression error can be written:
\begin{displaymath}
\sigma^2(x) = \phi(x) \left( \frac{\Phi^T \Phi}{n\tau} + \Lambda^{-1} \right)^{-1}\phi(x)^T
\end{displaymath}
Since $\bar{p}$ is finite, by the strong law of large numbers,  the $\bar{p} \times \bar{p}$ matrix $\frac{1}{n}\Phi^T\Phi$ converges almost surely as $n \rightarrow \infty$. We so have the following almost sure convergence:
\begin{equation}\label{degenerate_as}
\sigma^2(x)  \stackrel{n \rightarrow \infty}{\longrightarrow}  \sum_{p \leq \bar{p}} \frac{ \tau \lambda_p}{\tau + \lambda_p}  \phi_p(x)^2
\end{equation}
$\blacksquare$ 

\subsection{Proof of Theorem \ref{theo1}: the lower  bound for $\sigma^2(x)$}

The objective is to  find a lower bound for $\sigma^2(x)$ for non-degenerate kernels. Let us consider the Karhunen-Lo\`eve decomposition of  
$
Z(x) = \sum_{p \geq 0}Z_p \sqrt{\lambda_p} \phi_p(x)
$ where $(Z_p)_p$ is a sequence of independent Gaussian random variables with mean zero and variance 1.
If we denote by $a_i(x)$ the coefficients of the BLUP associated to $Z(x)$, the mean squared error can be written

\begin{eqnarray*}
\sigma^2(x) & = &  \mathbb{E}\left[ \left( Z(x) - \sum_{i=1}^{n} a_i(x) Z(x_i)\right)^2\right] \\
 & = & \mathbb{E}\left[ \left(  \sum_{p \geq 0}  \sqrt{\lambda_p} \left(  \phi_p(x) - \sum_{i=1}^{n} a_i(x) \phi_p(x_i) \right) Z_p  \right)^2\right] \\
 & = &  \sum_{p \geq 0} \lambda_p \left( \phi_p(x) - \sum_{i=1}^{n} a_i(x) \phi_p(x_i) \right)^2 \\
\end{eqnarray*}
Then, for a fixed $\bar{p}$, the following inequality holds:
\begin{equation}\label{ineq_sig_sigLUPA}
\sigma^2(x) \geq \sum_{p \leq \bar{p}} \lambda_p \left( \phi_p(x) - \sum_{i=1}^{n} a_i(x) \phi_p(x_i) \right)^2 = \sigma^2_{LUP,\bar{p}}(x)
\end{equation}
$\sigma^2_{LUP,\bar{p}}(x)$ is the MSE of the LUP of coefficients $a_i(x)$ associated to the Gaussian process $Z_{\bar{p}}(x) = \sum_{p \leq \bar{p}}Z_p \sqrt{\lambda_p} \phi_p(x)$. Let us consider $\sigma^2_{\bar{p}}(x)$ the MSE of the BLUP of  $Z_{\bar{p}}(x)$, we have  the following inequality: 
\begin{equation}\label{ineq_sigLUP_sigBLUPA}
\sigma^2_{LUP,\bar{p}}(x) \geq  \sigma^2_{\bar{p}}(x)
\end{equation}
Since $Z_{\bar{p}}(x)$ has a degenerate kernel,   the almost sure convergence given in equation (\ref{degenerate_as}) holds for $\sigma^2_{\bar{p}}(x)$. Then, considering inequalities (\ref{ineq_sig_sigLUPA}) and (\ref{ineq_sigLUP_sigBLUPA}) and the convergence (\ref{degenerate_as}), we obtain:
\begin{equation}\label{lowerboundA}
\liminf_{n \rightarrow \infty} \sigma^2(x) \geq \sum_{p \leq \bar{p}} \left( \frac{\tau \lambda_p}{\tau +  \lambda_p} \right) \phi_p(x)^2
\end{equation}
Taking the limit $\bar{p} \rightarrow \infty$ in the right hand side gives the desired result.
$\blacksquare$ 

\subsection{Proof of Theorem \ref{theo1}: the upper bound for $\sigma^2(x)$}

The objective is to find an upper bound for $\sigma^2(x)$. Since $\sigma^2(x)$ is the MSE of the BLUP associated to $Z(x)$, if we consider any other LUP associated to $Z(x)$ its MSE denoted by $\sigma_{LUP}^2(x)$ satisfies the following inequality:
\begin{equation}
\sigma^2(x) \leq \sigma_{LUP}^2(x)
\end{equation}
The idea is to find a LUP so that its MSE is a tight upper bound of $\sigma^2(x)$. Let us consider the LUP:
\begin{equation}\label{ALUP}
\hat{f}_{LUP}(x) = k(x)^T A z^{n}
\end{equation}
with $A$ the $n \times n$ matrix defined by 
$
A = L^{-1}+\sum_{k=1}^q(-1)^k(L^{-1}M)^kL^{-1}
$
with  $L = n \tau I + \sum_{p \leq p^*} \lambda_p [\phi_p(x_i) \phi_p(x_j)]_{1 \leq i,j \leq n}$, $M = \sum_{p > p^*} \lambda_p [\phi_p(x_i) \phi_p(x_j)]_{1 \leq i,j \leq n}$, $q$ a finite integer and $p^*$ such that $\lambda_{p^*} < \tau$. The matrix $A$ is an approximation of the inverse of the matrix $L+M = n \tau I + K $. Then, the  MSE of the LUP (\ref{ALUP}) is given by:
\begin{displaymath}
\sigma^2_{LUP}(x) = k(x,x) - k(x)^T\left( 2A - A(n \tau I + K)A\right)k(x)
\end{displaymath}
and by substituting  the expression of $A$ into  the previous equation we obtain:
\begin{equation}\label{sigLUPA}
\sigma^2_{LUP}(x) = k(x,x) - k(x)^T L^{-1} k(x) - \sum_{i=1}^{2q+1} (-1)^ik(x)^T(L^{-1}M)^iL^{-1}k(x)
\end{equation}

First, let us consider the term $k(x)^TL^{-1}k(x)$. Since $p^* < \infty$, the matrix $L$ can be written:
\begin{equation}
L    =   n\tau I + \Phi_{p^*} \Lambda \Phi_{p^*}^T
\end{equation}
where $\Lambda = \mathrm{diag}(\lambda_i)_{1 \leq i \leq p^*}$,   $ \Phi_{p^*} = \left(\begin{array}{ccc} \phi(x_1)^T  &  \dots  &  \phi(x_{n})^T \end{array} \right)^T$ and  $ \phi(x) = (\phi_1(x), \dots, \phi_{p^*}(x))$.

Thanks to the Woodbury-Sherman-Morrison formula, the matrix $L^{-1}$ is given by:
\begin{equation}\label{iL}
L^{-1} = \frac{I}{n\tau } - \frac{\Phi_{p^*} }{n\tau }\left( \frac{\Phi_{p^*} ^T \Phi_{p^*} }{n\tau } + \Lambda^{-1} \right)^{-1} \frac{\Phi_{p^*} ^T}{n\tau}
\end{equation}
From the continuity of the inverse operator for invertible $p^* \times p^*$ matrices and by applying the strong law of large numbers, we obtain the following almost sure convergence :
\begin{eqnarray*}
k(x)^T L^{-1}k(x)  & =  &   \frac{1}{n\tau }\sum_{i=1}^{n} k(x,x_i)^2 - \frac{1}{\tau^2} \sum_{p,q=0}^{p^*} \left[\left( \frac{\Phi_{p^*} ^T \Phi_{p^*} }{n\tau } + \Lambda^{-1} \right)^{-1} \right]_{p,q} \\
	& \times & \left[ \frac{1}{n} \sum_{i=1}^{n} k(x,x_i)\phi_p(x_i) \right]  \left[ \frac{1}{n} \sum_{j=1}^{n} k(x,x_j)\phi_q(x_j) \right]   \\
\end{eqnarray*}
\begin{displaymath}
 \stackrel{n \rightarrow \infty}{\longrightarrow}   \frac{1}{\tau} \mathbb{E}_\mu [k(x,X)^2] - \frac{1}{\tau^2} \sum_{p,q = 0}^{p^*}\left[\left( \frac{ I}{\tau } + \Lambda^{-1} \right)^{-1}\right]_{p,q} \mathbb{E}_\mu [k(x,X)\phi_p(X)] \mathbb{E}_\mu [k(x,X)\phi_q(X)] 
\end{displaymath}
where $\mathbb{E}_\mu$ is the expectation with respect to the design measure $\mu$.
We note that we can use the Woodbury-Sherman-Morrison formula and the strong law of large numbers since $p^*$ is finite and independent of $n$. 
Then, the orthonormal property of the basis $(\phi_p(x))_{p \geq 0}$ implies:
\begin{displaymath}
\mathbb{E}_\mu [k(x,X)^2] = \sum_{p \geq 0} \lambda_p^2 \phi_p(x)^2, \qquad  \mathbb{E}_\mu [k(x,X)\phi_p(X)] =   \lambda_p \phi_p(x)
\end{displaymath}
Therefore,  we have the following almost sure convergence:
\begin{equation}\label{kiLk}
k(x)^T L^{-1}k(x) \stackrel{n \rightarrow \infty}{\longrightarrow} \sum_{p \leq p^*} \frac{ \lambda_p^2}{ \lambda_p + \tau} \phi_p(x)^2 + \frac{1}{\tau} \sum_{p > p^*} \lambda_p^2 \phi_p(x)^2
\end{equation}

Second, let us consider the term $\sum_{i=1}^{2q+1} (-1)^ik(x)^T(L^{-1}M)^iL^{-1}k(x)$. We have the following equality:
\begin{eqnarray*}
k(x)^T(L^{-1}M)^iL^{-1}k(x)  & = &  \sum_{l=0}^i \left( \begin{array}{c} i \\ l \end{array} \right) \frac{1}{n \tau}  k(x)^T \left( \frac{M}{n \tau} \right)^l \left(-\frac{L'M}{(n\tau)^2} \right)^{i-l}k(x) \\
 & & -k(x)^T \left( \frac{M}{n \tau} \right)^l \left(-\frac{L'M}{(n\tau)^2} \right)^{i-l}\frac{L'}{(n \tau)^2}k(x) \\
\end{eqnarray*}
where:
\begin{equation}
L'   =  \Phi_{p^*} \left( \frac{\Phi_{p^*} ^T \Phi_{p^*} }{n\tau } + \Lambda^{-1} \right)^{-1} \Phi_{p^*} ^T = \sum_{p,p' \leq p^*} d^{(n)}_{p,p'}[ \phi_p(x_i) \phi_p(x_j)]_{1 \leq i,j \leq n}
\end{equation}
with $ d_{p,p'}^{(n)} = \left[ \left( \frac{\Phi_{p^*}^T \Phi_{p^*}}{n\tau} + \Lambda^{-1} \right)^{-1} \right]_{p,p'}$. Since $q < \infty$, we can obtain the convergence in probability of $\sum_{i=1}^{2q+1} (-1)^ik(x)^T(L^{-1}M)^iL^{-1}k(x)$ from the ones of:
\begin{equation}
k(x)^T\frac{1}{n }\left(\frac{M}{n } \right)^j \left(\frac{L'M}{n^2}\right)^{i-j} k(x)
\end{equation}
and:
\begin{equation}
k(x)^T\left(\frac{M}{n } \right)^j \left(\frac{L'M}{n^2}\right)^{i-j}\frac{L'}{n^2} k(x)
\end{equation}
with $i \leq 2q+1$ and $ j \leq i$.
Let us consider $k(x)^T\frac{1}{n }\left(\frac{M}{n } \right)^j \left(\frac{L'M}{n^2}\right)^{i-j} k(x)$ and $i > j$, we have:
\begin{equation}\label{kMLMkA}
k(x)^T\frac{1}{n }\left(\frac{M}{n } \right)^j \left(\frac{L'M}{n^2}\right)^{i-j} k(x)   =   	\sum_{\substack{p_1,\dots, p_{i-j} \leq p^* \\ p_1',\dots, p_{i-j}' \leq p^*}}
		d^{(n)}_{p_1,p_1'} \dots d^{(n)}_{p_{i-j},p_{i-j}'}
		\sum_{\substack{q_1,\dots, q_{i-j} > p^* \\ m_1,\dots, m_j >  p^*}}
		S^{(n)}_{q,m}  \\
\end{equation}
 with:
\begin{eqnarray*}
S^{(n)}_{q,m} & = &  	
	\left(  \frac{\sqrt{\lambda_{m_1}}}{n } \sum_{r=1}^{n} k(x,x_r) \phi_{m_1}(x_r) \right)  
	\left(   \frac{\sqrt{\lambda_{m_j}}}{n }  \sum_{r=1}^{n} \phi_{m_j}(x_r) \phi_{p_1'}(x_r) \right) \\
 & & 
	\times \left(  \frac{\lambda_{q_{i-j}}}{n } \sum_{r=1}^{n} k(x,x_r) \phi_{q_{i-j}}(x_r) 
	\sum_{r=1}^{n} \phi_{p_{i-j}}(x_r) \phi_{q_{i-j}}(x_r) \right) \\
 & &	
	\times \prod_{l=1}^{j-1} \frac{\sqrt{ \lambda_{m_l} \lambda_{m_{l+1}}}}{n } \sum_{r=1}^{n} \phi_{m_l}(x_r) \phi_{m_{l+1}}(x_r)  
	\prod_{l=1}^{i-j-1}  \frac{\lambda_{q_l}}{n } \sum_{r=1}^{n} \phi_{q_l}(x_r) \phi_{p_{l+1}}(x_r)  \sum_{r=1}^{n} \phi_{q_l}(x_r) \phi_{p_{l}'}(x_r) 
	\\
\end{eqnarray*}
We consider now the term:
\begin{equation}
a_{q,p,p'}^{(n)} =    \frac{ \lambda_{q}}{n } \sum_{r=1}^{n} \phi_{q}(x_r) \phi_{p}(x_r) 
	  \frac{1}{n } \sum_{r=1}^{n} \phi_{p'}(x_r) \phi_{q}(x_r)  
\end{equation}
with $p,p' \leq p^*$.
From  Cauchy Schwarz inequality and thanks to the following inequality:
\begin{displaymath}
|\phi_p(x)|^2  \leq  \frac{1}{\lambda_p} \sum_{p' \geq 0}\lambda_{p'}|\phi_{p'}(x)|^2  =  \lambda_p^{-1}k(x,x) 
\end{displaymath}
we obtain (using $\lambda_p \geq \lambda_{p^*}$, $\forall p \leq p^*$ and $[\sum_{r=1}^n|\phi_q(x_r)|]^2 \leq n\sum_{r=1}^n\phi_q(x_r)^2$):
\begin{displaymath}
\left|  a^{(n)}_{q,p,p'} \right|  \leq 
 \sigma^2 \lambda_{p^*}^{-1}  \frac{\lambda_q}{n} \sum_{r=1}^{n}\phi_q(x_r)^2  \qquad  \forall p,p' \leq  p^*
\end{displaymath}
with $\sigma^2 = \sup_{x} k(x,x)$. 
Considering the expectation with respect to the distribution of points  $x_r$, we obtain $\forall \bar{p} < \infty$:
\begin{displaymath}
\mathbb{E}_\mu \left[ \sum_{q > \bar{p}} \left|  a^{(n)}_{q,p,p'} \right| \right] \leq   \sigma^2 \lambda_{p^*}^{-1}   \sum_{q > \bar{p}} \lambda_q
\end{displaymath}
From  Markov inequality, $\forall \delta > 0$, we have:
\begin{equation}\label{Markov_ineq1}
\mathbb{P}_\mu \left( \left| \sum_{q > \bar{p}}a^{(n)}_{q,p,p'} \right| > \delta\right) \leq \frac{\mathbb{E}_\mu \left[\left|\sum_{q > \bar{p}}a^{(n)}_{q,p,p'} \right| \right]}{\delta} \leq \frac{ \sigma^2 \lambda_{p^*}^{-1}   \sum_{q > \bar{p}} \lambda_q}{\delta}
\end{equation}

Furthermore, $\forall \delta >0$, $\forall \bar{p} > p^*$: 
\begin{displaymath}
\mathbb{P}_\mu\left( \left| \sum_{q >  p^*}a^{(n)}_{q,p,p'} \right| > 2\delta\right) \leq 
\mathbb{P}_\mu\left( \left| \sum_{ p^* < q \leq \bar{p}}a^{(n)}_{q,p,p'} \right| > \delta \right)
+
\mathbb{P}_\mu\left( \left| \sum_{q > \bar{p}}a^{(n)}_{q,p,p'} \right| > \delta  \right)
\end{displaymath}
We have for all $q  \in (p^*, \bar{p}] : a^{(n)}_{q,p,p'} \rightarrow a_{q,p,p'} = \lambda_q    \delta_{q = p} \delta_{q = p'} = 0$ (with $\delta$ the Kronecker product), as $n \rightarrow \infty$,  therefore:
\begin{displaymath}
\limsup_{n \rightarrow \infty} 
\mathbb{P}_\mu \left( \left| \sum_{q > p*}a^{(n)}_{q,p,p'} \right| > 2 \delta\right) \leq \frac{ \sigma^2 \lambda_{p^*}^{-1}   \sum_{q > \bar{p}} \lambda_q}{\delta}
\end{displaymath}

Taking the limit $\bar{p} \rightarrow \infty$ in the right hand side, we obtain the  convergence in probability of  $\sum_{q > p^*}a^{(n)}_{q,p,p'} $ when $n \rightarrow \infty$:
\begin{equation}\label{probconv1}
  \sum_{q > p^*} \frac{ \lambda_{q}}{n } \sum_{r=1}^{n} \phi_{q}(x_r) \phi_{p}(x_r) 
	  \frac{1}{n } \sum_{r=1}^{n} \phi_{p'}(x_r) \phi_{q}(x_r)  
 \stackrel{\mathbb{P}_\mu}{\longrightarrow}
   0
\qquad \forall p,p' \leq p^*
\end{equation}
Following the same method, we obtain the  convergence:
\begin{equation}\label{probconv2}
\sum_{q >  p^*}  \frac{\lambda_{q}}{n } \sum_{r=1}^{n} k(x,x_r) \phi_{q}(x_r) 
	\sum_{r=1}^{n} \phi_{p}(x_r) \phi_{q}(x_r) 
 \stackrel{\mathbb{P}_\mu}{\longrightarrow}
 0
\qquad \forall p \leq p^*
\end{equation}
Let us return to $ S^{(n)}_{q,m} $. By using  Cauchy Schwarz inequality and  bounding by the constant $M$ all the terms independent of $q_i$ and $m_i$, we obtain:
\begin{eqnarray*}
\left| \sum_{\substack{q_1,\dots, q_{i-j}  >  p^*}} S^{(n)}_{q,m} \right| & \leq  &  
		M
		\prod_{l=1}^{j}  \lambda_{m_l}
		 \frac{1}{n } \sum_{r=1}^{n} \phi_{m_l}(x_r)^2 \\
 & & 
	\times \left| \sum_{\substack{q_{i-j} > p^*}} \left(  \frac{\lambda_{q_{i-j}}}{n } \sum_{r=1}^{n} k(x,x_r) \phi_{q_{i-j}}(x_r) 
	\sum_{r=1}^{n} \phi_{p_{i-j}}(x_r) \phi_{q_{i-j}}(x_r) \right) \right| \\
 & &	
	\times \left| \sum_{\substack{q_1,\dots, q_{i-j-1} > p^*}} \prod_{l=1}^{i-j-1}  \frac{\lambda_{q_l}}{n } \sum_{r=1}^{n} \phi_{q_l}(x_r) \phi_{p_{l+1}}(x_r)  \sum_{r=1}^{n} \phi_{q_l}(x_r) \phi_{p_{l}'}(x_r)  \right|
	\\
\end{eqnarray*}
Since $\sum_{p \geq 0} \lambda_p \phi_p(x)^2 = k(x,x) \leq \sigma^2$, we have the inequality $	0 \leq \sum_{m_1,\dots,m_{j}} \prod_{l=1}^{j}  \lambda_{m_l}
		 \frac{1}{n } \sum_{r=1}^{n} \phi_{m_l}(x_r)^2
		\leq ( \sigma^2)^j$. 
Thus, for $i > j$ and from (\ref{probconv1}) and (\ref{probconv2}) we obtain the following convergence in probability when $n \rightarrow \infty$:
\begin{eqnarray*}
		\sum_{\substack{q_1,\dots, q_{i-j} >  p^* \\ m_1,\dots, m_j >  p^*}}
		S^{(n)}_{q,m}
		 &  \stackrel{\mathbb{P}_\mu}{\longrightarrow} &  	0
		 \\
\end{eqnarray*}
Therefore, from (\ref{kMLMkA}) we obtain the following convergence when $n \rightarrow \infty$:
\begin{equation}\label{kMLMk0}
k(x)^T\frac{1}{n }\left(\frac{M}{n } \right)^j \left(\frac{L'M}{n^2}\right)^{i-j} k(x)  \stackrel{\mathbb{P}_\mu}{\longrightarrow} 	0
\qquad \forall i < j
\end{equation}
Following the same guideline as previously, it can be shown that when $n \rightarrow \infty$:
\begin{equation}\label{kMLMLk}
 k(x)^T\frac{1}{n }\left(\frac{M}{n } \right)^j \left(\frac{L'M}{n^2}\right)^{i-j} \frac{L'}{n^2} k(x) 
  \stackrel{\mathbb{P}_\mu}{\longrightarrow}   0 \qquad \forall i \leq j
\end{equation}

From the convergences (\ref{kMLMk0}) and (\ref{kMLMLk}), we deduce the following one when $n \rightarrow \infty$:
\begin{equation}\label{kLMqLk}
 k(x)^T\left(L^{-1}M\right)^qL^{-1} k(x)
  -  \frac{1}{n} k(x)^T\left( \frac{M}{n}\right)^q k(x) \stackrel{\mathbb{P}_\mu}{\longrightarrow} 0
\end{equation}
Therefore, to complete the proof we have to show that:
\begin{displaymath}
 \frac{1}{n} k(x)^T\left( \frac{M}{n}\right)^q k(x)   \stackrel{\mathbb{P}_\mu}{\longrightarrow}     \sum_{p >  p^*}\lambda_p^{q+2} \phi_p(x)^2
\end{displaymath}

Let us consider for a fixed $j \geq 1$:
\begin{displaymath}
 \frac{1}{n} k(x)^T\left( \frac{M}{n}\right)^j k(x)  =  \sum_{\substack{ m_1,\dots, m_j > p^*}} a^{(n)}_m(x)
\end{displaymath}
with $m = (m_1,\dots,m_j)$ and:
\begin{eqnarray*}
a^{(n)}_{m}(x) & = &
\left( \frac{1}{n} \sum_{r=1}^{n}k(x,x_r)\phi_{m_1}(x_r)\right)
\left( \frac{1}{n} \sum_{r=1}^{n}k(x,x_r)\phi_{m_j}(x_r)\right)
\\
 & & \times \prod_{l=1}^{j-1} \frac{1}{n } \sum_{r=1}^{n} \phi_{m_l}(x_r) \phi_{m_{l+1}}(x_r)  
	\prod_{i=1}^j \lambda_{m_i}
\\
\end{eqnarray*}

From  Cauchy-Schwarz inequality, we have:
\begin{eqnarray}\label{eq52} 
\left|a^{(n)}_{m}(x) \right| & \leq &
\left( \frac{1}{n} \sum_{r=1}^{n} k(x,x_r)^2 \right) \prod_{i=1}^j \frac{1}{n} \sum_{r=1}^{n} \lambda_{m_i} \phi_{m_i}(x_r)^2 \\ \label{eq53}
 & \leq &  \sigma^4 \prod_{i=1}^j \frac{1}{n} \sum_{r=1}^{n} \lambda_{m_i} \phi_{m_i}(x_r)^2 
\end{eqnarray}

Therefore, considering the expectation with respect to  the distribution of the points $(x_r)_{r=1,\dots,n}$, we have:
\begin{eqnarray*}
\mathbb{E}_{\mu} \left[\left|a^{(n)}_{m}(x) \right| \right] & \leq  &  \sigma^4
\left( \prod_{i=1}^j  \lambda_{m_i} \right) \frac{1}{n^j} \sum_{t_1,\dots,t_j=1}^{n}  \mathbb{E}_{\mu} \left[ \phi_{m_1}(X_{t_1})^2 \dots \phi_{m_j}(X_{t_j})^2  \right]
\qquad \forall x \in \mathbb{R}^d
\end{eqnarray*}
The following inequality holds  uniformly in $t_1,\dots,t_j = 1,\dots, n$:
\begin{eqnarray*}
\mathbb{E}_{\mu} \left[ \prod_{i=1}^j  \phi_{m_i}(X_{t_i})^2  \right]  & \leq & b_{m}
\\
\end{eqnarray*}
where $b_{m} = \sum_{ \substack{  \mathcal{P} \in \Pi( \{1, \dots, j\}) \\  \mathcal{P} = \cup_{r=1}^l I_r}}
	\prod_{r=1}^l \mathbb{E}_{\mu}\left[ \prod_{i \in I_r} \phi_{m_i}(X)^2 \right] $ 
because the term  of left hand side of the inequality  is equal to one of the terms in the sum of the right hand side. Here
 $\Pi( \{1, \dots, j\}) $  is the collection of all partitions of $\{1, \dots, j\}$ and $I_r \cap  I_{r'} =  \emptyset$, $\forall r \neq r'$. We hence have:
\begin{eqnarray*}
\mathbb{E}_{\mu} \left[\left|a^{(n)}_{m}(x)\right| \right] & \leq  & \sigma^4  
\prod_{i=1}^j  \lambda_{m_i} b_{m}
\end{eqnarray*}
Since $\sum_{p \geq 0} \lambda_p \phi_p(x)^2 \leq \sigma^2 $, we have:
\begin{eqnarray*}
\sum_{m_1,\dots,m_j > p^*}
\prod_{i=1}^j  \lambda_{m_i}  b_{m} &  =  &   \sum_{m_1,\dots,m_j > p^*}
	\prod_{l=1}^j  \lambda_{m_l}
	\sum_{ \substack{  \mathcal{P} \in \Pi( \{1, \dots, j\}) \\  \mathcal{P} = \cup_{r=1}^l I_r}}
	\prod_{r=1}^l \mathbb{E}_{\mu}\left[ \prod_{i \in I_r} \phi_{m_i}(X)^2 \right]  \\
 &  =  &  
	\sum_{ \substack{  \mathcal{P} \in \Pi( \{1, \dots, j\}) \\  \mathcal{P} = \cup_{r=1}^l I_r}}
	\prod_{r=1}^l \mathbb{E}_{\mu}\left[ \prod_{i \in I_r} \sum_{m_i  > p^*} \lambda_{m_i} \phi_{m_i}(X)^2 \right]  \\
& \leq & \sigma^{2j} \# \{  \Pi( \{1, \dots, j\}) \} \\
\end{eqnarray*}
Since the cardinality of the collection $\Pi( \{1, \dots, j\})$ of  partitions of  $\{1, \dots, j\}$ is finite,  the series  $\sum_{m_1,\dots,m_j > p^*}
\prod_{i=1}^j  \lambda_{m_i} b_{m}$ converges.
Furthermore, as it is a series with non-negative terms, $\forall \varepsilon > 0$, $\exists \bar{p} > p^*$ such that :
\begin{displaymath}
   \sigma^4 \sum_{\substack{ m  \in M^C_{\bar{p}} }} \prod_{i=1}^j \lambda_{m_i}b_m \leq \varepsilon
\end{displaymath}
where $M^C_{\bar{p}}$ designs the complement of $M_{\bar{p}}$ defined by the collection of  $ m = (m_1,\dots,m_j)$ such that:
\begin{displaymath}
M = \{ m = (m_1,\dots,m_j) \quad \mathrm{such \, that} \quad m_i > p^*, \quad i = 1,\dots,j\}
\end{displaymath}
\begin{displaymath}
M_{\bar{p}} = \{ m = (m_1,\dots,m_j) \quad \mathrm{such \, that} \quad   p^* < m_i  \leq \bar{p}, \quad i = 1,\dots,j\}
\end{displaymath}
\begin{displaymath}
M^C_{\bar{p}} = M \setminus M_{\bar{p}}
\end{displaymath}

Therefore, we have $\forall \delta > 0$, $\forall \varepsilon >0$ $\exists \bar{p} > 0$ such that uniformly in $n$:
\begin{displaymath}
  \sum_{\substack{ m  \in M^C_{\bar{p}}  }} \mathbb{E}_{\mu} \left[\left|a^{(n)}_{m}(x)\right| \right]  \leq  \frac{\varepsilon \delta}{2}
\end{displaymath}
Applying  the Markov inequality, we obtain:
\begin{equation}\label{anMC}
   \mathbb{P} \left( \sum_{\substack{ m \in M^C_{\bar{p}} }} \left|a^{(n)}_{m}(x)\right| > \frac{\delta}{2}  \right)  \leq \varepsilon
\end{equation}
Furthermore, by denoting $a_{m}(x) = \lim_{n \rightarrow \infty} a^{(n)}_{m}(x)$, we have:
\begin{equation}\label{limam}
  a_{m}(x) =  \lambda_{m_1} \lambda_{m_j} \phi_{m_1}(x)\phi_{m_j}(x)  \prod_{i=1}^{j}\lambda_{m_i} \prod_{i=1}^{j-1}\delta_{m_i= m_{i+1}}
\end{equation}
and from  Cauchy-Schwarz inequality (see equation (\ref{eq53})), we have:
\begin{displaymath}
 | a_{m}(x)| \leq \sigma^4   \prod_{i=1}^j \lambda_{m_i}
\end{displaymath}
We hence can deduce the inequality:
\begin{equation}\label{aMC}
     \sum_{\substack{ m  \in M^C_{\bar{p}} }}\left|   a_{m}(x)  \right|  \leq   \sigma^4  \sum_{\substack{ m  \in M^C_{\bar{p}} }}   \prod_{i=1}^j \lambda_{m_i} 
\end{equation}
Thus, $\exists \bar{p}$ such that  $  \sum_{\substack{ m \in M^C_{\bar{p}} }}  \left|    a_{m}(x) \right|   \leq  \frac{\delta}{2}$ for all $x \in  \mathbb{R}^d$.
From the inequalities  (\ref{anMC}) and  (\ref{aMC}), we find that $\exists \bar{p}$ such that:
\begin{displaymath}
  \mathbb{P}_\mu \left( \left|  \sum_{\substack{ m \in M }}a^{(n)}_{m}(x) -   \sum_{\substack{ m \in M}} a_{m}(x) \right|  > 2\delta \right) \leq  \varepsilon + \mathbb{P}_\mu \left( \left|  \sum_{\substack{ m \in M_{\bar{p}} }}a^{(n)}_{m}(x) -   \sum_{\substack{ m \in M_{\bar{p}} }} a_{m}(x) \right|  > \delta \right)
\end{displaymath}
Since $M_{\bar{p}}$ is a finite set:
\begin{displaymath}
  \limsup_{n \rightarrow \infty}   \mathbb{P}_\mu \left( \left|  \sum_{\substack{ m \in M_{\bar{p}} }}a^{(n)}_{m}(x) -   \sum_{\substack{ m \in M_{\bar{p}} }} a_{m}(x) \right|  > \delta \right) = 0
\end{displaymath}
therefore:
\begin{displaymath}
  \limsup_{n \rightarrow \infty}  \mathbb{P}_\mu \left( \left|  \sum_{\substack{ m \in M }}a^{(n)}_{m}(x) -   \sum_{\substack{ m \in M }} a_{m}(x) \right|  > 2 \delta \right) \leq  \varepsilon   
\end{displaymath}

The previous inequality holds $\forall \varepsilon > 0$,   thus we have the convergence in probability of $ \sum_{\substack{ m \in M }}a^{(n)}_{m}(x)$ to $   \sum_{\substack{ m \in M}} a_{m}(x)$ with (by using the limit in the equation (\ref{limam})):
\begin{displaymath}
 \sum_{\substack{ m \in M }} a_{m}(x) =  \sum_{p > p^*} \lambda_p^{j+2} \phi_p(x)^2
\end{displaymath}
Finally, we have the following convergence in probability when $n \rightarrow \infty$:
\begin{equation}
k(x)^T (L^{-1}M)^iL^{-1}k(x) \stackrel{n \rightarrow \infty}{\longrightarrow} \left( \frac{1}{\tau}\right)^{i+1} \sum_{p > p^*} \lambda_p^{i+2}\phi_p(x)^2
\end{equation}
We highlight that we cannot use the strong law of large numbers here due to the infinite sum in $M$.

 From the equation (\ref{sigLUPA}) and  the convergences (\ref{kiLk}) and (\ref{kLMqLk}), we obtain the following convergence in probability:
\begin{equation}
\sigma^2_{LUP}(x)  \stackrel{n \rightarrow \infty}{\longrightarrow}  \sum_{p \geq 0} \left( \lambda_p - \frac{ \lambda_p^2}{\tau +  \lambda_p} \right) \phi_p(x)^2 - \sum_{p > p^*} \lambda_p^2 \frac{\left( \frac{ \lambda_p}{\tau} \right)^{2q+1}}{\tau +  \lambda_p} \phi_p(x)^2
\end{equation}
By considering the limit $q \rightarrow \infty$ and the inequality $ \lambda_{p^*} < \tau$, we obtain the following upper bound for $\sigma^2(x)$:
\begin{equation}\label{upperboundA}
\limsup_{n \rightarrow \infty} \sigma^2(x) \leq \sum_{p \geq 0}  \frac{\tau \lambda_p}{\tau +  \lambda_p}   \phi_p(x)^2
\end{equation}
$\blacksquare$

\bibliographystyle{apalike}
\bibliography{biblio}

\label{fin}

\end{document}